
\input amstex
\documentstyle{amsppt}
\mag=1200
\magnification=\magstep1
\hcorrection{.25in}
\advance\vsize-.75in
\input epsf.tex

\def\Fig#1#2#3{\medskip\vbox{\epsfysize=#1truein \centerline{\epsfbox{#2}}
	\centerline{\bf Figure #3}}\medskip}

\def\hatP{\widehat P}
\def\hatQ{\widehat Q}
\def\rint{\text{\rm int }}
\def\nhd{\text{\rm nhd}}
\def\bp{\text{\bf p}}
\def\bq{\text{\bf q}}
\def\ep{\varepsilon}
\def\real{{\Bbb R}}
\def\que{{\Bbb Q}}
\def\zed{{\Bbb Z}}
\def\B{{\Cal B}}
\def\C{{\Cal C}}
\def\D{{\Cal D}}
\def\W{{\Cal W}}
\def\chix{\raise.5ex\hbox{$\chi$}}
\def\scoprod{\ \raise.5ex\hbox{$\scriptstyle\coprod$}\ }

\topmatter
\title Combinatorial methods in Dehn surgery\endtitle
\author C. McA. Gordon\endauthor
\leftheadtext{C. M{\fiverm c}A. GORDON}
\address Department of Mathematics, The University of Texas at Austin,
Austin, TX 78712-1082
\endaddress
\curraddr 
MSRI, 1000 Centennial Dr., Berkeley, CA 94720
\endcurraddr 
\email gordon\@math.utexas.edu\endemail
\thanks Partially supported by NSF Grant DMS-9626550.\endthanks
\thanks Research at MSRI supported in part by NSF grant \#DMS 9022140. 
\endthanks
\thanks To appear in the Proceedings of the Workshop of {\it Knots
'96}, the fifth MSJ International Research Institute in Knot Theory, 
Waseda University, Tokyo, Japan, 1996.
\endthanks

\abstract
This is an expository paper, in which we give a summary of some of the
joint work of John Luecke and the author on Dehn surgery.  We consider
the situation where we have two Dehn fillings $M(\alpha)$ and
$M(\beta)$ on a given 3-manifold $M$, each containing a surface that
is either essential or a Heegaard surface.  We show how a
combinatorial analysis of the graphs of intersection of the two
corresponding punctured surfaces in $M$ enables one to find faces of
these graphs which give useful topological information about
$M(\alpha)$ and $M(\beta)$, and hence, in certain cases, good upper
bounds on the intersection number $\Delta(\alpha, \beta)$ of the
two filling slopes.
\endabstract

\endtopmatter

\document 

\head 1. Introduction\endhead 

This expository paper is an expanded version of part of the author's
lectures at the Workshop.  The remaining part of the lectures covered
material that is discussed in Sections~2 and 3 of the survey article
\cite{Go}.

Here, we give an introductory account of the combinatorial methods
developed by John Luecke and the author (mainly in the proof of the
Knot Complement Conjecture \cite{GL2}) for studying various questions
about Dehn surgery.  For more details, see [GL2], [GL3] and [GL4].

The setting is that we have two Dehn fillings $M(\alpha)$, $M(\beta)$ on 
a 3-manifold $M$, which contain surfaces $\hatP,\hatQ$ respectively, and 
the main goal is to obtain sharp upper bounds on the intersection number 
of $\alpha$ and $\beta$, under certain hypotheses on $\hatP$ and $\hatQ$. 

In Section 2 we describe how $\hatP$ and $\hatQ$ give rise to labelled 
intersection graphs $G_P \subset \hatP$, $G_Q\subset \hatQ$. 
The idea of using such graphs in the context of Dehn surgery is due to 
Litherland \cite{L}. 

In Section 3 we introduce the algebraic notion of a {\it type\/}, which 
plays a key role in the theory. 

In Section 4 we show how certain faces or collections of faces of $G_P$ 
give topological information about $M(\beta)$ (and/or $\hatQ$), and of 
course similarly for $G_Q$. 

Section 5 contains the main combinatorial assertion, Theorem 5.1, which 
says that, if $G_P$ and $G_Q$ are graphs as in Section~2, then either 
$G_P$ represents all types or $G_Q$ contains a configuration we call a 
{\it $(p-\chix (\hatP))$-web\/}. 

Finally, Section 6 contains some  applications, mainly to Dehn surgery on 
knots in $S^3$, although we also say something about the case where 
$M(\alpha)$ and $M(\beta)$ are reducible. 

We have tried to keep the discussion as general as possible, allowing 
$\hatP$ and $\hatQ$ to have arbitrary genus until the last part of         
Section~5. 
It must be admitted that this generality is somewhat bogus, inasmuch as in 
all the applications in Section~6, $\hatQ$ is a 2-sphere and $\hatP$ is 
either a 2-sphere or a torus. 
(For statements which apply to surfaces of arbitrary genus see \cite{T} and 
\cite{R}.) 
The main point is that if $\hatQ$ is a 2-sphere then, if $G_Q$ contains a 
$(p-\chix(\hatP))$-web, it contains one that is innermost in the obvious 
sense, and in such a {\it great\/} web one can find useful faces, namely 
Scharlemann cycles. 
Nevertheless, we hope that presenting the material in this way clarifies 
the logical structure of the argument. 

\head 2. Graphs of surface intersections\endhead 

Let $M$ be a compact, orientable, irreducible 3-manifold with torus boundary.  
Let $\alpha$ be a {\it slope\/}, that is, the isotopy class of an essential 
embedded circle, on $\partial M$. 
Then the result of {\it $\alpha$-Dehn filling\/} on $M$ is the manifold 
$M(\alpha) = M\cup V_\alpha$, where $V_\alpha$ is a solid torus, glued to 
$M$ via a homeomorphism $\partial M\to \partial V_\alpha$ taking $\alpha$ 
to the boundary of a meridian disk of $V_\alpha$. 

If $\alpha$ and $\beta$ are two slopes on $\partial M$, then $\Delta(\alpha,
\beta)$ will denote the minimal geometric intersection number of $\alpha$ 
and $\beta$. 

One finds in practice that statements of the following form tend to hold, 
and this is the main type of result that we shall be seeking: 
\smallskip
\itemitem{} {\it if $M(\alpha)$ and $M(\beta)$ have certain special 
properties, then $\Delta (\alpha,\beta)$ is bounded above by some 
$\Delta_0$.} 
\smallskip

\noindent  
The bound $\Delta_0$ will depend on the properties in question. 
The situation in which we are particularly interested, and to which our 
methods apply, is when $M(\alpha)$ and $M(\beta)$ contain certain kinds 
of surfaces. 

So suppose that $\hatP ,\hatQ$ are (closed, orientable) surfaces 
in $M(\alpha)$, $M(\beta)$ respectively. 
We may assume that $\hatP$ meets $V_\alpha$ in a finite collection of 
meridian disks, so that $P = \hatP\cap M$ is a surface in $M$ each of 
whose boundary components has slope $\alpha$. 
Similarly, $\hatQ$ gives rise to a surface $Q$ in $M$ whose boundary 
components have slope $\beta$. 

For our machinery to be applicable, it is important that $\partial P$ and 
$\partial Q$ be non-empty. 
So, we may state our basic assumption as follows: 
\roster
\item"(A1)" {\it there are properly embedded surfaces $(P,\partial P)$, 
$(Q,\partial Q) \subset (M,\partial M)$ such that $\partial P$ and 
$\partial Q$ are non-empty, and each component of $\partial P$ 
(resp. $\partial Q$) has slope $\alpha$ (resp. $\beta$).}
\endroster

By an isotopy of (say) $P$ we may assume in addition that 
\roster
\item"(A2)" {\it $P$ and $Q$ intersect transversely, and each component of 
$\partial P$ meets each component of $\partial Q$ in exactly $\Delta (\alpha,
\beta)$ points.} 
\endroster

It is clear, however, that for any $\alpha$ and $\beta$ we can always find 
surfaces $P$ and $Q$, of arbitrary genus, satisfying (A1) and (A2), and 
therefore we must impose some additional conditions if we are to extract any 
non-trivial information from this situation. 
The technical assumption which enables us to do this is: 
\roster
\item"(A3)" {\it each arc component of $P\cap Q$ is essential in $P$ and 
in $Q$.}
\endroster

By (A2), $P\cap Q$ consists of a finite disjoint union of circles and 
properly embedded arcs, the endpoints of the arcs being the points of 
intersection of $\partial P$ with $\partial Q$. 
Condition (A3) says that for no arc component $\gamma$ of $P\cap Q$ is there 
an arc $\delta \subset \partial P$, with $\partial \gamma=\partial \delta$, 
and a disk $D\subset P$ such that $\partial D= \gamma\cup\delta$ (and 
similarly for $Q$). 

We now describe some conditions under which (A3) holds. 

Recall that a closed orientable 
surface $S$ in a 3-manifold  $N$ is {\it essential\/} 
if either $S$ has positive genus, is incompressible in $N$, and is not 
parallel to a component of $\partial N$, or $S$ is a 2-sphere which does 
not bound a 3-ball in $N$. 

First suppose 
\roster
\item"(a)" {\it $\hatP$ is essential in $M(\alpha)$ and $M$ contains no 
essential surface homeomorphic to $\hatP$.} 
\endroster

Let $K_\alpha \subset M(\alpha)$ denote the core of the solid torus $V_\alpha$. 
We choose $\hatP$ (among all essential surfaces in $M(\alpha)$ 
homeomorphic to $\hatP$) so that $p= |\partial P| = |\hatP\cap K_\alpha|$ 
is minimal. 
Then $p>0$ by hypothesis. 
Also, standard arguments show that $P$ is incompressible and boundary 
incompressible in $M$. 
In particular, we may assume that no circle component of $P\cap Q$ bounds a 
disk in $Q$. 

If $\hatQ$ also satisfies (a) (with $M(\alpha)$ replaced by $M(\beta)$), and 
is chosen to minimize $q= |\partial Q| = |\hatQ \cap K_\beta|$, then,  
again by standard arguments, we may assume that condition (A3) holds. 

Another case of interest is 
\roster
\item"(b)" {\it $\hatP$ is a Heegaard surface for $M(\alpha)$ and $K_\alpha$ 
cannot be isotoped to lie on $\hatP$.}
\endroster

If $\hatP$ satisfies (a) and $\hatQ$ satisfies (b) (for $M(\beta)$), then 
Gabai \cite{Ga2} shows that (A3) can be achieved if $K_\beta$ is put in 
{\it thin position\/} with respect to $\hatQ$. 
(Gabai explicitly treats the case when $M(\beta) \cong S^3$ and $\hatQ \cong 
S^2$, but his argument carries over verbatim to any Heegaard surface $\hatQ$.) 

Finally, if $\hatP$ and $\hatQ$ both satisfy (b) (for $M(\alpha)$ and 
$M(\beta)$ respectively), then Rieck has shown \cite{R}, using thin position, 
that again (A3) can be assumed to hold. 
(The case $M(\alpha) \cong M(\beta) \cong S^3$, $\hatP,\hatQ\cong S^2$ is 
done in \cite{GL2}, and was also proved independently by Gabai (unpublished).) 

Given surfaces $P,Q$ in $M$ satisfying (A1), (A2) and (A3), we focus on the 
arc components of $P\cap Q$, as they lie in $P$ and in $Q$. 
We regard these arcs as defining graphs $G_P$ and $G_Q$ in $\hatP$ and $\hatQ$ 
respectively, in the obvious way. 
Thus the (fat) vertices of $G_P$ are the disks $\hatP -\rint P$, the edges 
of $G_P$ are the arc components of $P\cap Q$ as they lie in $P$, and 
similarly for $G_Q$. 

We now encode some additional structure, as follows. 
Number the components of $\partial P$\ \ $1,2,\ldots,p$ in the order in 
which they appear on $\partial M$, and similarly number the components of 
$\partial Q$\ \ $1,2,\ldots,q$. 
This gives a corresponding numbering of the vertices of $G_P$ and $G_Q$. 
Consider an arc component of $P\cap Q$. 
Its endpoints are points of $\partial P\cap \partial Q$, the intersections 
of (say) the $i$th component of $\partial P$ with the $j$th component of 
$\partial Q$, and the $i'$th component of $\partial P$ with the $j'$th 
component of $\partial Q$. 
We then label the endpoints of the corresponding edge of $G_P$ with $j$ and 
$j'$ at vertices $i$ and $i'$ respectively, and similarly for $G_Q$. 
See Figure~1. 

\bigskip

\Fig{.75}{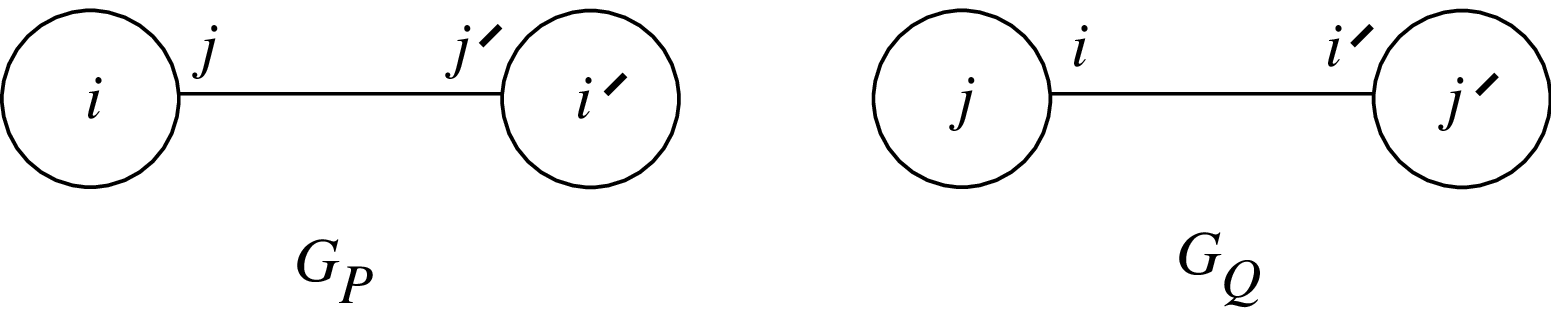}{1}

We next give a {\it sign\/}, $+$ or $-$, to each vertex of $G_P$, according 
to the direction on $\partial M$ of the orientation of the corresponding 
component of $\partial P$, induced by some chosen orientation of $P$. 
Equivalently, this is the sign of the corresponding intersection point 
of $K_\alpha$ with $\hatP$ (with respect to some chosen orientations of 
$M,\hatP$ and $K_\alpha$). 
In particular, if $\hatP$ separates $M(\alpha)$ then the signs of the vertices 
$1,2,\ldots,p$ of $G_P$ alternate. 
Similarly, we give a sign to each vertex of $G_Q$. 
Note that around each vertex of $G_P$ the edge-endpoint labels occur in order 
$1,2,\ldots,q$, $1,2,\ldots,q,\ldots$, repeated $\Delta=\Delta(\alpha,\beta)$  
times, the ordering being (say) anticlockwise at a positive vertex and 
clockwise at a negative vertex. 
See Figure~2.
Similarly, around each vertex of $G_Q$ we see the labels $1,2,\ldots,p$, 
$1,2,\ldots,p,\ldots,$ repeated $\Delta$ times. 

\bigskip

\Fig{2.5}{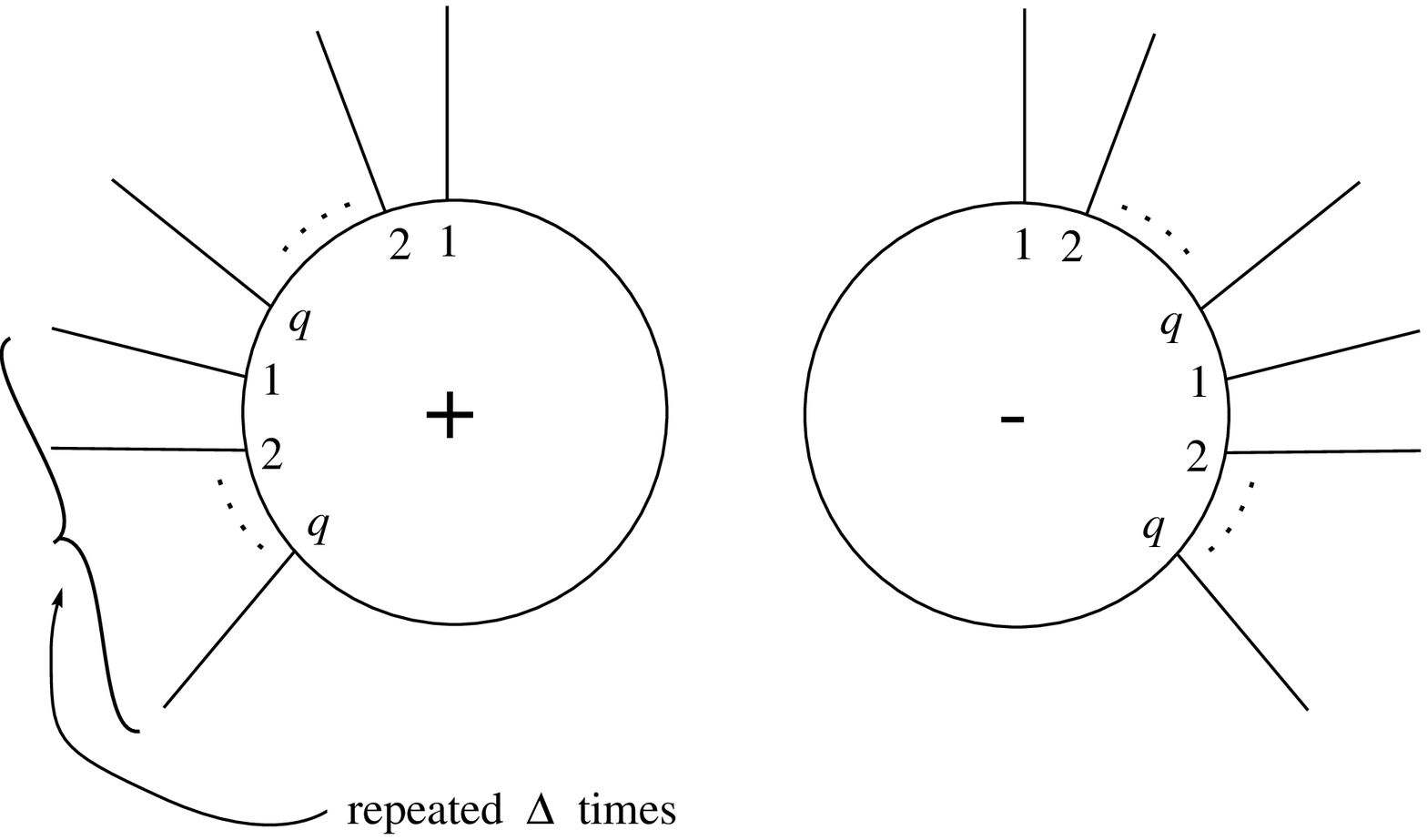}{2}

Since $M,P$ and $Q$ are orientable, an arc component of $P\cap Q$ joins 
points of intersection of $\partial P$ with $\partial Q$ of opposite sign. 
Hence we have the {\it parity rule\/}: if an edge of $G_P$ joins vertices 
$i$ and $i'$ and the corresponding edge of $G_Q$ joins  vertices $j$ and $j'$, 
then $i$ and $i'$ have the same sign if and only if $j$ and $j'$ have 
opposite signs. 
(See Figure~1.) 

Condition (A3) becomes, in graph-theoretic language, the condition that 
$G_P$ and $G_Q$ contain no trivial loops. 

We will denote by $\bq$ the set of (edge-endpoint) labels $\{1,2,\ldots,q\}$ 
of $G_P$. 
We also have the associated {\it $\bq$-intervals\/} $(1,2),(2,3),\ldots,
(q-1,q),(q,1)$. 

\bigskip

\Fig{2.5}{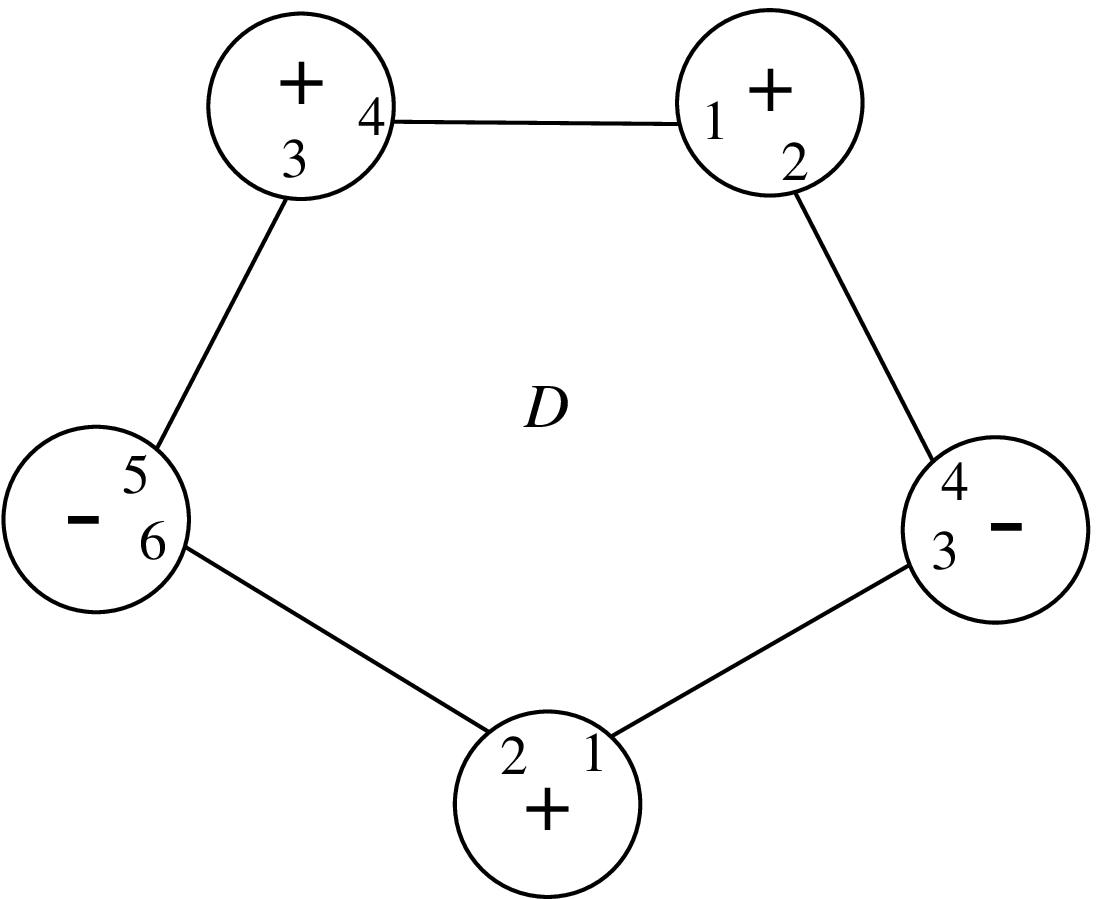}{3}

Note that if $D$ is a disk face of $G_P$, then $\partial D$ consists of 
an alternating sequence of {\it edges\/} and {\it corners\/}, where the 
edges are edges of $G_P$ (i.e., arc components of $P\cap Q$), and the 
corners are $\bq$-intervals $(i,i+1)$ (i.e., arcs in $\partial P$ between 
consecutively labelled components $i,i+1$ of $\partial Q$. 
See Figure~3. 

Similar remarks apply to $G_Q$. 

We may describe the philosophy behind our approach as consisting of two 
parts: 
\medskip 

\noindent {\bf Combinatorial:}
use the one-one correspondence between the edges of the  labeled graphs 
$G_P$ and $G_Q$ to show that either $G_P$ contains certain configurations 
or $G_Q$ contains certain (perhaps other) configurations.
\smallskip

\noindent{\bf Topological:}
use the existence of certain collections of faces of $G_P$ (resp. $G_Q$) 
to get topological information about the pair $(M(\beta),\hatQ)$ 
(resp. $(M(\alpha),\hatP))$.

These will be elaborated in Sections 5 and 4 respectively. 

\remark{Remark} 
For simplicity we have assumed that $\partial M$ is a torus, but everything 
carries over without much change if we allow $M$ to have additional 
boundary components (with $\alpha$ and $\beta$ being slopes on some fixed 
torus component of $\partial M$). 
In this more general setting, we may also allow $\hatP$ and/or $\hatQ$ 
to have non-empty boundary.
\endremark 

\head 3. Types\endhead

We now introduce a purely algebraic concept, that of a {\it type\/}, 
which plays an important role in the theory under discussion. 
The reader may choose to skip to Section~4 and return to the present 
section when necessary. 

Let $q$ be a positive integer. 
A {\it $q$-type\/} is a $q$-tuple ($\ep_1,\ep_2,\ldots,\ep_q)$, where 
$\ep_i=\pm$, $1\le i\le q$. 

We say that $a= (a_1,a_2,\ldots,a_q) \in \zed^q$ {\it represents\/} the 
$q$-type $\tau = (\ep_1,\ep_2,\ldots,\ep_q)$ if and only if 
\roster
\item"(1)" for some $\eta=\pm$, sign $a_i = \eta\ep_i$ for each $i$ such 
that $a_i\ne0$; and 
\item"(2)" $\sum\limits_{i=1}^q |a_i|\ge2$.
\endroster

Condition (2) merely excludes the zero vector, together with vectors all 
of whose co-ordinates are zero except one, which is $\pm1$. 

Note that $a$ represents $\tau$ if and only if it represents $-\tau$. 
Also, there are $2^{q-1}$ $q$-types up to sign, including the {\it trivial\/} 
type, which by definition has $\ep_1=\ep_2 =\cdots = \ep_q$. 

\example{Example} 
There are $2^2=4$ 3-types (up to sign): $+++$, $-++$, $+-+$, and $++-$. 
$(1,0,-2)$ represents $-++$ and $++-$. 
\endexample

A set $A\subset \zed^q$ {\it represents all $q$-types\/} if and only if, 
for each $q$-type $\tau$, there exists $a\in A$ such that $a$ represents 
$\tau$. 

\example{Examples} 
\roster
\item"(1)" If $|k|\ge2$, then $\{(0,\ldots,0,k,0,\ldots,0)\}$ represents 
all $q$-types.
\item"(2)" $\{(1,0,1),(-1,1,0),(1,1,-1)\}$ represents all 3-types. 
\item"(3)" If $1\le r\le q$, then $A\subset\zed^q$ represents all $q$-types 
if and only if 
\item"{}" $\{(a_1,\ldots,a_r,0,a_{r+1},\ldots,a_q): 
(a_1,\ldots,a_q)\in A\}$ represents all $(q+1)$-types.
\endroster
\endexample 

For $A\subset \zed^q$, let $c(A)$ be the number of co-ordinates $1\le i\le q$ 
such that $a_i\ne0$ for some $a\in A$. 
Then we may regard $A\subset \zed^{c(A)} \subset \zed^q$ in the obvious way. 
By (3) above, $A$ represents all $c(A)$-types if and only if it represents 
all $q$-types. 
Thus we may use the phrase $A$ {\it represents all types\/} without 
ambiguity. 

The key property of types is the following ``all types implies torsion'' 
theorem due to Parry \cite{P}. 

\proclaim{Theorem 3.1 {\rm (Parry)}} 
If $A\subset \zed^q$ represents all types then there exists $A_0 \subset A$ 
such that $\zed^q$ modulo the subgroup generated by $A_0$ has non-trivial 
torsion. 
\endproclaim 

We remark that it is not true that $A_0$ can be chosen to represent all types; 
equivalently, it is not true that if $A\subset \zed^q$ is a minimal set 
of representatives of all types then $\zed^q/(A)$ has non-trivial torsion. 
(Here, {\it minimal\/} means that no proper subset of $A$ represents all 
types.) 
For example, $A= \{(1,2,1), (-1,1,1), (1,-1,2), (1,1,-1)\}$ 
is clearly a minimal set of representatives of all types, but 
$\zed^3 /(A)=0$. 

A useful fact, proved in \cite{GL4, Lemma 4.4}, is that if $A\subset \zed^q$ 
represents all types then either $A$ contains a basis for $\real^q$ or there 
exists $A_1\subset A$ such that $c(A_1) <q$ and $A_1$ represents all types. 
As an immediate consequence we have 

\proclaim{Lemma 3.2} 
If $A$ is a minimal set of representatives of all types then $A$ contains 
a basis for $\real^{c(A)}$.
\endproclaim 

\head 4. The topology\endhead 

Consider the solid torus $V_\beta \subset M(\beta)$. 
The surface $\hatQ$ meets $V_\beta$ in $q$ meridian disks, cutting $V_\beta$ 
into $q$ 3-balls. 
For each $\bq$-interval $\lambda = (i,i+1)$, let $H_\lambda$ denote the 
3-ball (1-handle) consisting of that part of $V_\beta$ between the meridian 
disks whose boundaries are components $i$ and $i+1$ of $\partial Q$. 

Let $D$ be a disk face of $G_P$. 
Then $D\subset M\subset M(\beta)$, where the edges of $\partial D$ are arcs 
in $Q$, and the corners $(i,i+1)$ of $\partial D$ are arcs in $\partial P$ 
running along the corresponding 1-handle $H_{(i,i+1)}$. 
See Figure~4, which shows $\partial D$ for the face $D$ illustrated 
in Figure~3. 

\bigskip

\Fig{2.5}{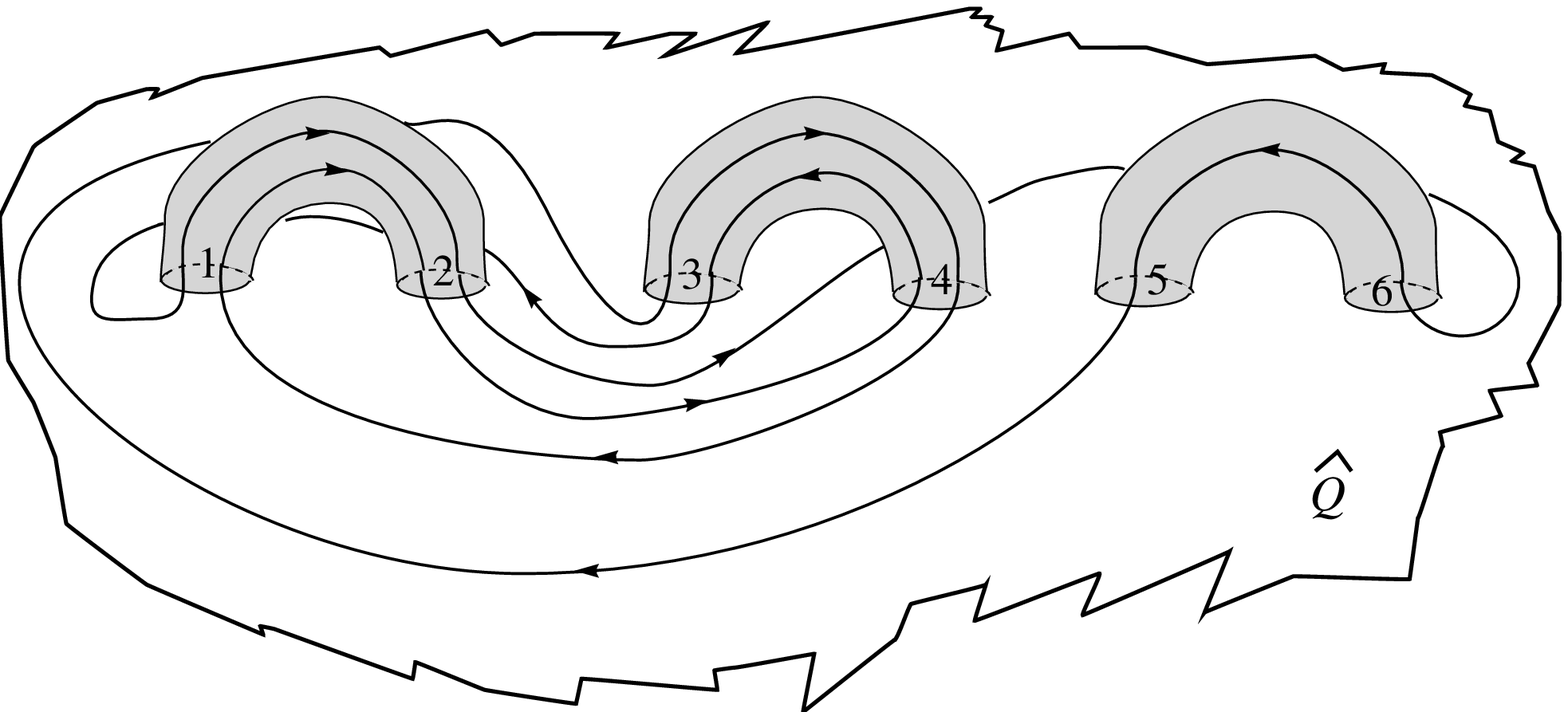}{4}

Let $\D$ be a set of disk faces of $G_P$. 
Let $c(\D)$ be the set of $\bq$-intervals $\lambda = (i,i+1)$ that occur as 
corners of faces in $\D$. 
Let $e(\D)$ be the set of edges of $G_Q$ corresponding to edges belonging 
to faces in $\D$. 

Suppose that the edges $e(\D)$ are contained in a disk $E\subset \hatQ$. 
(Note that this will necessarily be the case if $\hatQ$ is a 2-sphere.) 
We take $E$ to also contain the fat vertices of $G_Q$ at the endpoints 
of the edges in $e(\D)$. 

Recall that if $\hatQ$ is an essential surface, then we may assume that 
no circle component of $P\cap Q$ bounds a disk in $Q$. 
Hence $(\text{int }\D)\cap  \hatQ =\emptyset$. 
If $\hatQ$ is a Heegaard surface, this need not hold. 
Assume for the moment, however, that $(\text{int }\D)\cap E=\emptyset$; 
we shall see later (at the end of the present section) that in the cases 
of interest this can be achieved after some disk exchanges on int~$\D$. 

Define a 3-manifold $N(E,\D)\subset M(\beta)$ by 
$$N(E,\D) = \nhd \biggl( E\cup \bigcup_{\lambda\in c(\D)} H_\lambda \cup \D
\biggr)\ .$$ 

Note that $N(E,\D)$ has a natural handle decomposition, with a single 0-handle 
(corresponding to $E$), the 1-handles $H_\lambda$, and 2-handles whose cores 
are the elements of $\D$. 
Also, if $W$ denotes the handlebody $\nhd (E\cup \bigcup_{\lambda\in c(\D)} 
H_\lambda)$, then we can read off the element of $\pi_1(W)$ (and hence 
of $H_1(W)$) represented by the boundary of a disk $D\in \D$ by reading the 
sequence of corners of $D$ (with sign, given by the sign of the corresponding 
vertex of $G_P$) as we go around $\partial D$ in some direction. 

The idea is to find $\D$ such that $N(E,\D)$, or perhaps the surface 
$\partial N(E,\D)$, gives topological information about $M(\beta)$, or 
the surface $\hatQ$, possibly a contradiction. 

The simplest example of this philosophy is when $\D$ consists of a single 
face of $G_P$ which is a {\it Scharlemann cycle\/}. 
This is a disk face $D$ of $G_P$ such that all the vertices of 
$D$ have the  same 
sign, and all the corners of $D$ are the same $\bq$-interval 
$\lambda = (i,i+1)$. 
See Figure~5. 
Scharlemann cycles first appeared (in a slightly different context) in 
\cite{S1}. 
The following theorem, which is contained in \cite{S2, Proof of 
Proposition~5.6}, is immediate; see Figure~6. 
(Note that vertices $i$ and $i+1$ of $G_Q$ are of opposite sign, by the 
parity rule.) 

\bigskip

\Fig{1.5}{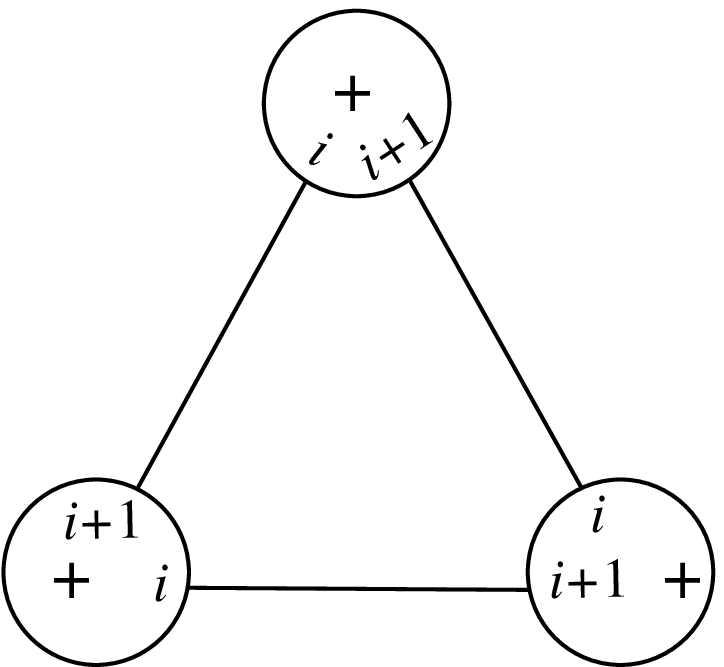}{5}

\proclaim{Theorem 4.1} 
Let $D$ be a Scharlemann cycle in $G_P$, such that the edges $e(D)$ lie in 
a disk $E\subset \hatQ$. 
Then $N(E,D)$ is a punctured lens space.
\endproclaim 

Note that $\pi_1(N(E,D))\cong \zed_k$, where $k$ is the {\it length\/} 
of $D$, i.e., the number of edges in $\partial D$. 

\bigskip

\Fig{2.0}{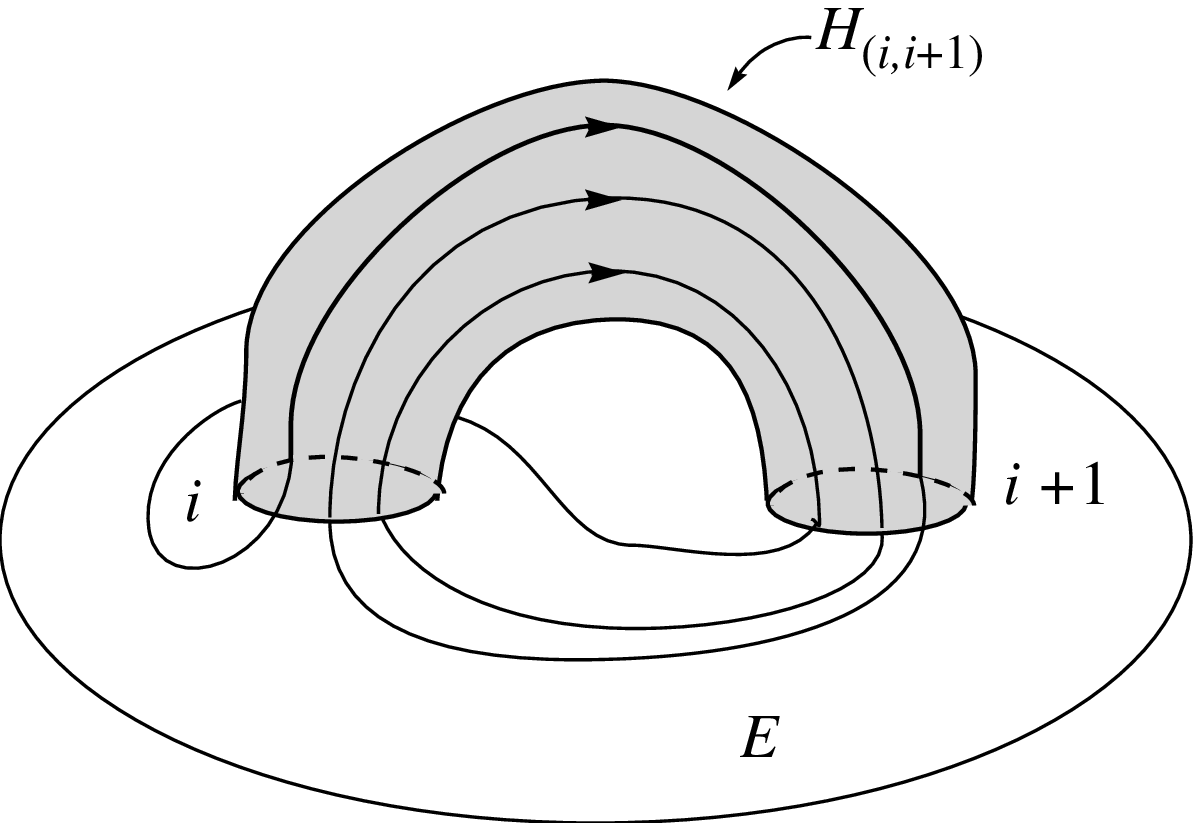}{6}

Let us examine the implications for $(M(\beta),\hatQ)$ of the existence of 
a Scharlemann cycle in $G_P$ as in Theorem~4.1. 
\roster
\item"(Sa)" If $\hatQ$ is a Heegaard surface, we conclude tht $M(\beta)$ has 
a lens space summand.
\endroster

Now suppose that $\hatQ$ is an essential surface. 
Note that (after a small isotopy of $N(E,D))$ we may assume that 
$\partial N(E,D)\cap \hatQ =E$.
Let $E'$ be the disk $\overline{\partial N(E,D)-E}$, and let 
$\hatQ' = (\hatQ-E)\cup E'$.
Then $\hatQ' \cong\hatQ$. 
Also, $|\hatQ' \cap K_\beta| = q-2$, since the two points of intersection 
of $K_\beta$ with $\hatQ$ corresponding to the 
labels $i$ and $i+1$ that appear in the corners of the Scharlemann cycle 
$D$, have been eliminated. 

We therefore have the following. 
\roster
\item"(Sb)" If $\hatQ$ is an incompressible surface of positive genus, then 
$\hatQ'$ is also incompressible. 
But this is a contradiction, because of the minimality of $q$. 
\item"(Sc)" If $\hatQ$ is an essential 2-sphere, then $M(\beta)$ has a lens 
space summand. 
In fact, by the minimality of $q$, $\hatQ$ decomposes $M(\beta)$ as 
$M'\# $ (lens space).
\item"(Sd)" For a variant of (Sc), let us say that a 2-sphere  $S$ in a 
3-manifold $N$ is {\it $\que$-essential\/} if $S$ does not bound a 
$\que$-homology ball in $N$. 
(In particular, this is the case if $S$ is non-separating.) 
It is clear that if $\hatQ$ is $\que$-essential in $M(\beta)$, then so is 
$\hatQ'$. 
Therefore, if we choose $\hatQ$ so that $q$ is minimal over all 
$\que$-essential 2-spheres in $M(\beta)$, then the existence of a 
Scharlemann cycle in $G_P$ is a contradiction. 
\endroster

Although Scharlemann cycles are very useful, their existence is not 
guaranteed in all situations of interest. 
A more global concept, which gives similar topological conclusions, 
is that a set of faces representing all types. 

A {\it $\bq$-type\/} is a $q$-type $\tau = (\ep_1,\ep_2,\ldots,\ep_q)$ 
(see Section~3), where the $i$th co-ordinate of $\tau$ is formally 
associated with the $i$th $\bq$-interval $(i,i+1)$. 

Let $D$ be a disk face of $G_P$. 
By taking the algebraic sum of the corners of $D$ we get an element 
$[D]\in\zed^q$, the free abelian group on the set $\bq$-intervals. 

We say that $D$ {\it represents\/} the $\bq$-type $\tau$ if and only if 
\roster
\item"(1)" $[D]$ represents $\tau$ (in the sense of Section 3); and 
\item"(2)" ({\it coherence condition\/}) no $\bq$-interval $(i,i+1)$ appears 
as a corner of $D$ at two vertices of opposite sign.
\endroster

As an example, the disk $D$ shown in Figure~3 fails to satisfy the coherence 
condition. 

A set of (disk) faces $\D$ of $G_P$ {\it represents all types\/} if 
and only if for 
each $\bq$-type $\tau$ there exists $D\in \D$ such that $D$ represents $\tau$. 
If such a set of faces exists we say that $G_P$ {\it represents all types\/}. 

\remark{Remark} 
The role of the coherence condition (2) is twofold. 
First, it ensures that condition (2) in the algebraic definition of 
representing a type (see Section~3) is automatically satisfied 
(since $G_P$ has no trivial loops). 
Second, it makes it possible to search for representatives of types in a 
more organized fashion, as we shall see in Section~5 (cf.\ Proposition~5.2). 
\endremark 

If $D$ is a face of $G_P$ such that small neighborhoods in $D$ of 
the edges $e(D)$ all lie on the same side of $\hatQ$, we say that $D$ 
{\it lies locally\/} on that side of $\hatQ$. 
To get the topological conclusions that we want from a set of 
representatives of all types, we first need to observe the following 
(cf.\ \cite{GL2, Proof of Proposition 3.2, part~(A)}). 

\proclaim{Lemma 4.2} 
Let $\D$ be a set of faces of $G_P$ representing all types. 
Then there exists $\D' \subset\D$ such that $\D'$ represents all types 
and the elements of $\D'$ all lie locally on the same side of $\hatQ$.
\endproclaim 

\demo{Proof} 
Label the two local sides of $\hatQ$ in $M(\beta)$, $B$ and $W$ (black and 
white). 

Let $\lambda  = (i,i+1)$ be a $\bq$-interval, and recall the definition 
of the 1-handle $H_\lambda$ at the beginning of this section. 
Orient (the core of) $H_\lambda$ from the end corresponding to component 
$i$ of $\partial Q$ to the end corresponding to component $i+1$ of 
$\partial Q$. 
Then we have the partition 
$$\{\text{$\bq$-intervals}\} = \B \scoprod \W\scoprod \B\W\scoprod \W\B\ ,$$ 
where 
$$\lambda \in  \left\{ \matrix \B\cr \W\cr \B\W\cr \W\B\endmatrix\right\} 
\text{ if $H_\lambda$ runs from }
\left\{ \matrix 
B\text{ to } B\cr 
W\text{ to } W\cr 
B\text{ to } W\cr 
W\text{ to } B 
\endmatrix\right\}\ .$$

Suppose the conclusion of the lemma is false. 
Then, if $\B \ne \emptyset$, there exists a $\B$-type $\tau_B$ such that no 
member of $\D$ represents $\tau_B$. 
Similarly, if $\W\ne \emptyset$, there exists a $\W$-type $\tau_W$ such that 
no member of $\D$ represents $\tau_W$. 
Define a $\bq$-type $\tau$ by 
$$\align
&\tau|\B = \tau_B\ , \ \text{ if } \B\ne\emptyset\quad ;\quad 
\tau|\W = \tau_W\ ,\ \text{ if } \W\ne \emptyset\ ;\cr 
&\tau|\lambda = +\ ,\ \forall\ \lambda\in \B\W\quad ;\quad 
\tau|\lambda = -\ , \ \forall\ \lambda \in \W\B\ .
\endalign$$ 

By hypothesis, there exists $D\in \D$ such that $D$ represents $\tau$. 
Note that, for some orientation of $\partial D$, every corner of $D$ 
in $\B\W\scoprod \W\B$ runs from $B$ to $W$, and hence $D$ can contain 
no such corner. 
Therefore either all corners of $D$ belong to $\B$, or all corners belong 
to $\W$. 
It follows that either $\B\ne\emptyset$ and $D$ represents $\tau_B$, or 
$\W\ne\emptyset$ and $D$ represents $\tau_W$, contradicting the definitions 
of $\tau_B$ and $\tau_W$.\qed
\enddemo

\proclaim{Theorem 4.3} 
Let $\D$ be a set of faces of $G_P$ representing all types, such that 
the edges $e(\D)$ lie in a disk $E\subset\hatQ$. 
Then there exists $\D' \subset \D$ such that the elements of $\D'$ 
all lie locally on the same side of $\hatQ$, 
$H_1(N(E,\D'))$ is finite 
and non-zero, and $\partial N(E,\D')\cong S^2$. 
\endproclaim 

\demo{Proof} 
(cf.\ \cite{GL4, proof of Proposition 4.5}). 
We may assume that $\D$ is minimal. 
Hence, by Lemma~4.2, the elements of $\D$ all lie locally on the same 
side of $\hatQ$. 
Also, by Lemma~3.2, there exists $\D'\subset \D$ such that $\{[D]:D\in\D'\}$ 
is a basis for $\real^{c(\D)}$. 
Let $N' = N(E,\D')$. 
Then it follows that $H_1(N')$ is finite and $\partial N'\cong S^2$. 

It remains to show that $H_1(N')\ne0$. 

Let $N= N(E,\D)= \nhd (N'\cup (\D-\D'))$. 
Since $\partial N'\cong S^2$, attaching to $N'$ the 2-handles whose cores 
are the elements of $\D-\D'$ is equivalent to removing from $N'$ an equal 
number of open 3-balls. 
By Theorem~3.1 there exists $\D_0 \subset\D$ such that $H_1(N_0)$ has 
non-trivial torsion, where $N_0 = N(E,\D_0)$. 
Note that $N_0\subset N\subset N'$. 
Suppose $H_1(N') =0$. 
Then $N'$ is a homology 3-ball, and the Mayer-Vietoris exact sequence 
gives $H_1(\partial N_0)\cong H_1(N_0) \oplus H_1(N'- N_0)$, 
contradicting the fact that $H_1(N_0)$ has non-trivial torsion. 
Hence $H_1(N')\ne0$, as desired.\qed
\enddemo 

\remark{Remark} 
Since  $H_1(N)\cong H_1(N')$, the above proof shows that if $\D$ is a 
minimal set of representatives of all types then $H_1(N(E,\D))$ is finite 
and non-zero. 
As was noted in Section~3, the analogous statement is not true in the 
purely algebraic setting. 
\endremark 

It follows from Theorem 4.3 that if $G_P$ represents all types, and the edges 
of the relevant faces are contained in a disk in $\hatQ$, then we get 
conclusions almost identical to (Sa), (Sb), (Sc) and (Sd) above for 
Scharlemann cycles.
\roster
\item"(Ta)" If $\hatQ$ is a Heegaard surface, then $M(\beta)$ has a 
summand $\hat N'$ with $H_1 (\hat N')$ finite and non-zero.
\item"(Tb)" If $\hatQ$ is an incompressible surface of positive genus, 
this is a contradiction.
(We get a new surface $\hatQ'$, obtained from $\hatQ$ by a disk exchange, 
with $|\hatQ'\cap K_\beta| = q- |c(\D)|$.) 
\item"(Tc)" If $\hatQ$ is an essential 2-sphere, then $\hatQ$ decomposes 
$M(\beta)$ as $M'\,\#\, \widehat N'$, where $H_1(\widehat N')$ 
is finite and non-zero. 
\item"(Td)" If $\hatQ$ is a $\que$-essential 2-sphere, this is a 
contradiction. 
\endroster

Finally, we return to the assumption, when $\hatQ$, is a Heegaard surface, 
that $(\text{int }\D)\cap E=\emptyset$.

Let $\D$ be a collection of faces of $G_P$ representing all types (this 
includes the case where $\D$ consists of a single Scharlemann cycle). 
We may assume that $(\text{int }\D)\cap E$ consists of a finite number 
of disjoint circles and arcs, properly embedded in $E$. 
Let $\gamma$ be such a circle or arc, and let $E_1$ and $E_2$ be the 
components of $E-\gamma$. 
There is a corresponding partition 
$$c(\D) = \C_1 \scoprod \C_2 \scoprod \C_{12}\ ,$$ 
where $\lambda \in \C_i$ if the 1-handle $H_\lambda$ has both ends in 
$E_i$, $i=1,2$, and $\lambda \in \C_{12}$ if $H_\lambda$ has one end in 
$E_1$ and one in $E_2$. 

Then, exactly as in the proof of Lemma~4.2 (the $c(\D)$-type $\tau$ here 
will be defined by orienting all $\lambda \in \C_{12}$ so as to 
run from (say) $E_1$ to $E_2$), there exists $\D' \subset \D$ such that 
$\D'$ represents all types and $c(\D')\subset \C_1$ (say). 
Since we may assume that $\D$ is a minimal set of representatives of all 
types, we conclude that the vertices of $G_Q$ corresponding to the labels that 
appear in the corners $c(\D)$ all lie in $E_1$. 

If $\gamma$ is an arc, or if $\gamma$ is a circle and $\bar E_1$ is a disk, 
we replace $E$ by $\bar E_1$ (moved slightly so as to eliminate the 
intersection $\gamma$ with $\text{int }\D$). 
If $\gamma$ is a circle and $\bar E_2$ is a disk, we use a standard 
innermost circle argument to cut-and-paste $\D$, using subdisks of 
$\bar E_2$, so as to reduce $(\text{int }\D)\cap E$.  
Thus we eventually get a collection of disks $\widetilde{\D}$, whose 
boundaries are the same as those in $\D$, such that 
$(\text{int }\widetilde{\D})\cap E=\emptyset$. 

\remark{Remark} 
There are situations when it is useful to consider submanifolds of 
$M(\beta)$ of the form $\nhd(E\cup \bigcup H_\lambda \cup \D)$ where $E$ 
is a subsurface of $\hatQ$ other than a disk. 
For example, the case where $\D$ consists of a Scharlemann 
cycle and $E$ is an annulus arises in \cite{GL3, Section~3}. 
\endremark 

\head 5. The combinatorics\endhead 

To state the main result of this section, Theorem~5.1 below, we need the 
following definition. 

A subgraph $\Lambda$ of $G_Q$ is a {\it $k$-web\/} ($k$ a non-negative 
integer) if and only if 
\roster
\item"(1)" all the vertices of $\Lambda$ have the same sign, and 
\item"(2)" at most $k$ edge-endpoints of $G_Q$ at vertices of $\Lambda$ 
are not endpoints of edges of $\Lambda$.
\endroster

See Figure 7 for an example. 

\bigskip

\Fig{3.0}{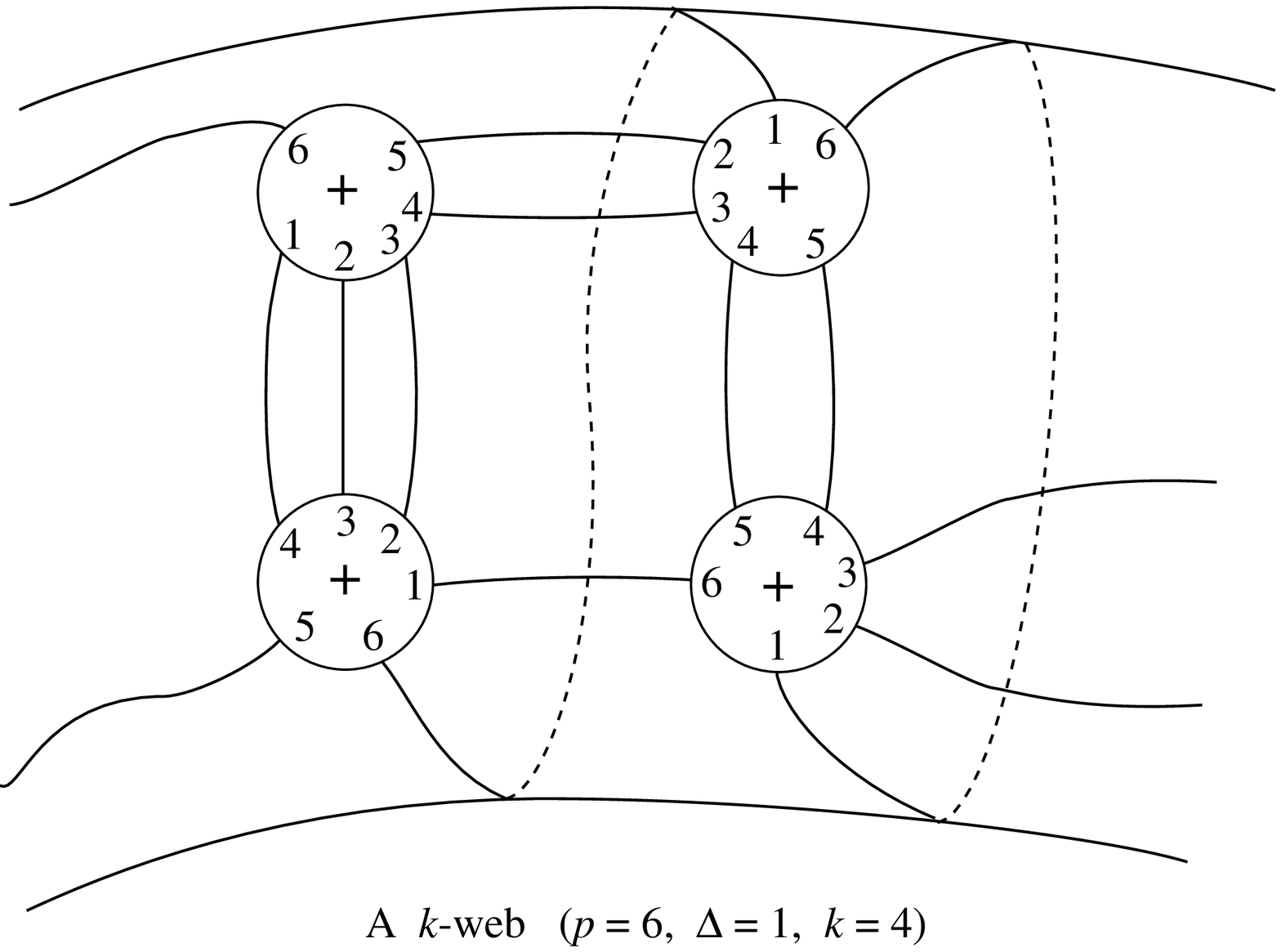}{7}

In the following theorem, $G_P$ and $G_Q$ are graphs as in Section~2. 

\proclaim{Theorem 5.1} 
Assume that $\Delta >1-\chix (\hatP)/p$. 
Then either $G_P$ represents all types or $G_Q$ contains a  
$(p-\chix (\hatP))$-web.
\endproclaim 

Most of the remainder of this section will be devoted to a sketch of 
the proof of Theorem~5.1. 

First we deal with the case of the trivial type, which is different from 
the general case. 

Define $\Lambda$ to be the subgraph of $G_Q$ consisting of all edges of 
$G_Q$ joining vertices of the same sign. 
If some component of $\Lambda$ is not a $(p-\chix (\hatP))$-web, then 
$G_Q$ contains more than $(p-\chix (\hatP))$ edges joining vertices of 
opposite sign. 
By the parity rule, the corresponding edges of $G_P$ join vertices of the 
same sign. 
Let $\Sigma$ be the subgraph of $G_P$ consisting of all vertices of 
$G_P$, and all edges of $G_P$ joining vertices of the same sign. 
Then 
$$\align 
\chix (\hatP) & = V-E + {\textstyle\sum}
\chix (f) \ \text{ (summed over all faces $f$ 
of $\Sigma$ )}\\ 
& < p- (p-\chix (\hatP)) + {\textstyle\sum} \chix (f)\ .
\endalign$$ 
Hence $\sum\chix (f) >0$, implying that $\Sigma$ has a disk face. 
This face is then a face of $G_P$ whose vertices all have the same sign, 
and which therefore represents the trivial $\bq$-type. 

For non-trivial $\bq$-types we proceed as follows. 

Define a graph $\Gamma \subset \hatP$ by: 

$\bullet$ the {\it vertices\/} of $\Gamma$ consist of the {\it fat\/} vertices 
(i.e., the vertices of $G_P$), together with {\it dual\/} vertices 
$v(D)$, one in the interior of each face $D$ of $G_P$; 

$\bullet$ the {\it edges\/} of $\Gamma$ join each dual vertex $v(D)$ to the fat 
vertices in the boundary of the corresponding face $D$. 

See Figure 8. 

\bigskip

\Fig{1.5}{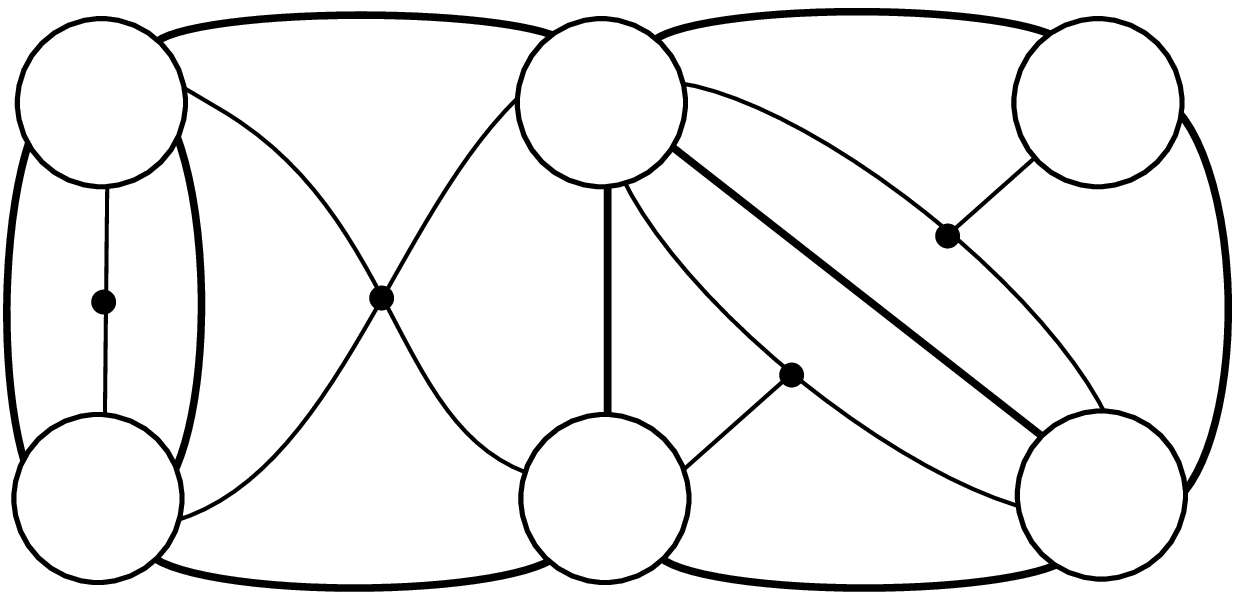}{8}

Now let $\tau$ be a $\bq$-type. 
We define the directed graph $\Gamma (\tau) \subset \hatP$ to be the graph 
$\Gamma$, with edges oriented according to the following rule. 
Let $e$ be an edge of $\Gamma$, with one endpoint at the fat vertex $v$, 
lying in a $\bq$-interval $(i,i+1)$ at that vertex. 
Then we orient $e$ 
\roster
\item"{}" {\it inwards\/} at $v$ if $(\tau|(i,i+1))\cdot\text{sign }v=+$, and
\item"{}" {\it outwards\/} at $v$ if $(\tau|(i,i+1))\cdot\text{sign }v=-$.
\endroster

As an example, if $q=6$ and $\tau = +-+- - -$, then around the vertices 
of $G_P$ the edges of $\Gamma(\tau)$ are oriented as shown in Figure~9.

\bigskip

\Fig{2.00}{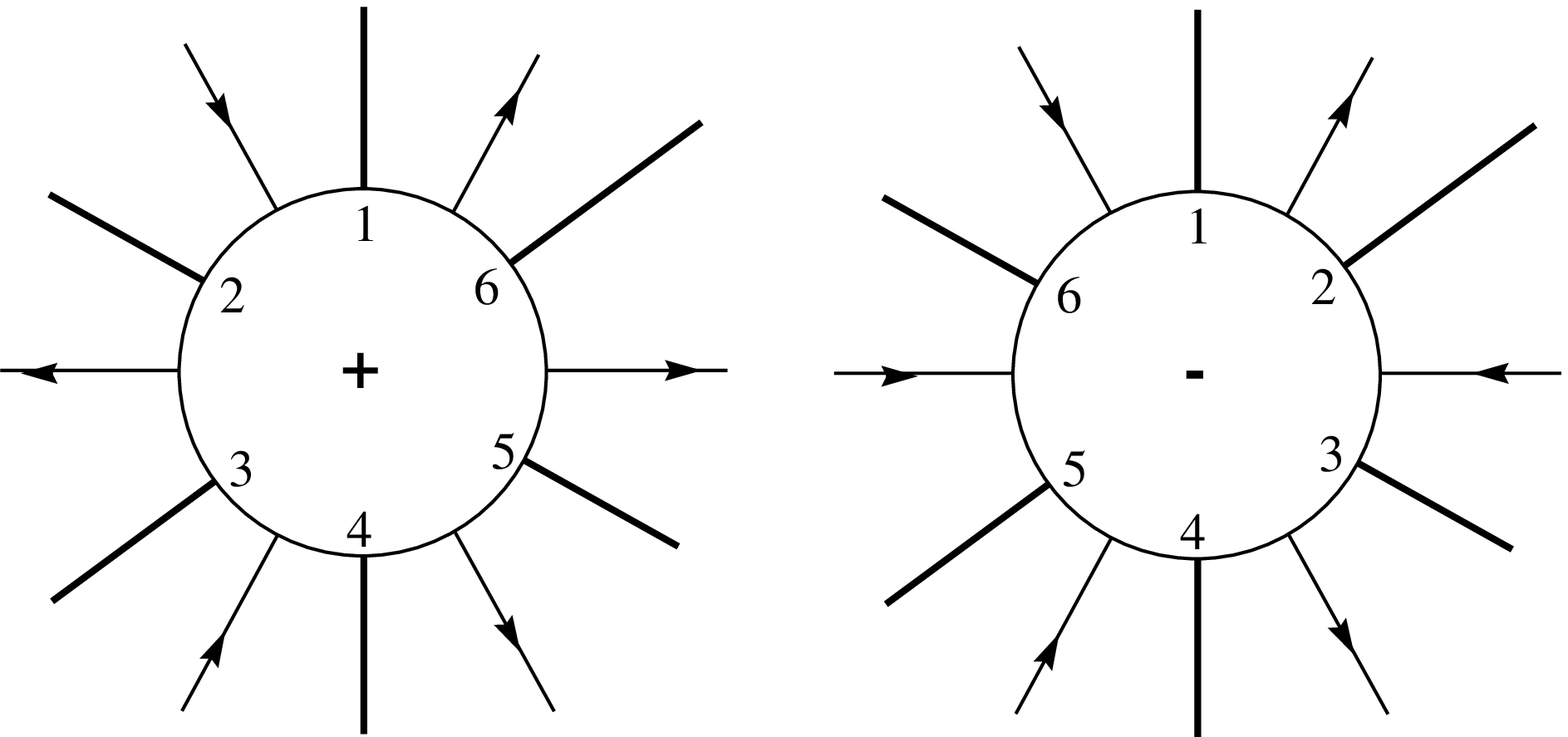}{9}

The reason for defining the orientation on the edges of $\Gamma (\tau)$ in 
this way is the following, which is immediate from the definition. 

\proclaim{Proposition 5.2} 
A disk face $D$ of $G_P$ represents $\tau$ if and only if the corresponding 
dual vertex $v(D)$ is a sink or source of $\Gamma (\tau)$.
\endproclaim 

Our search for representatives of types is thus translated into the study 
of the directed graphs $\Gamma (\tau)\subset \hatP$. 
Here we will use the following very elegant combinatorial formulation of 
the Poincar\'e-Hopf Index Theorem, due to Glass \cite{Gl}. 

Consider any directed graph $\Omega$ in a closed surface $S$. 
For each vertex $v$ of $\Omega$, let $s(v)$ be the number of {\it switches\/} 
(i.e., changes in orientation of successive edges) around $v$, and for 
each face $f$ of $\Omega$, let $s(f)$ be the number of switches around 
$\partial f$. 
See Figure~10. 

\bigskip

\Fig{1.5}{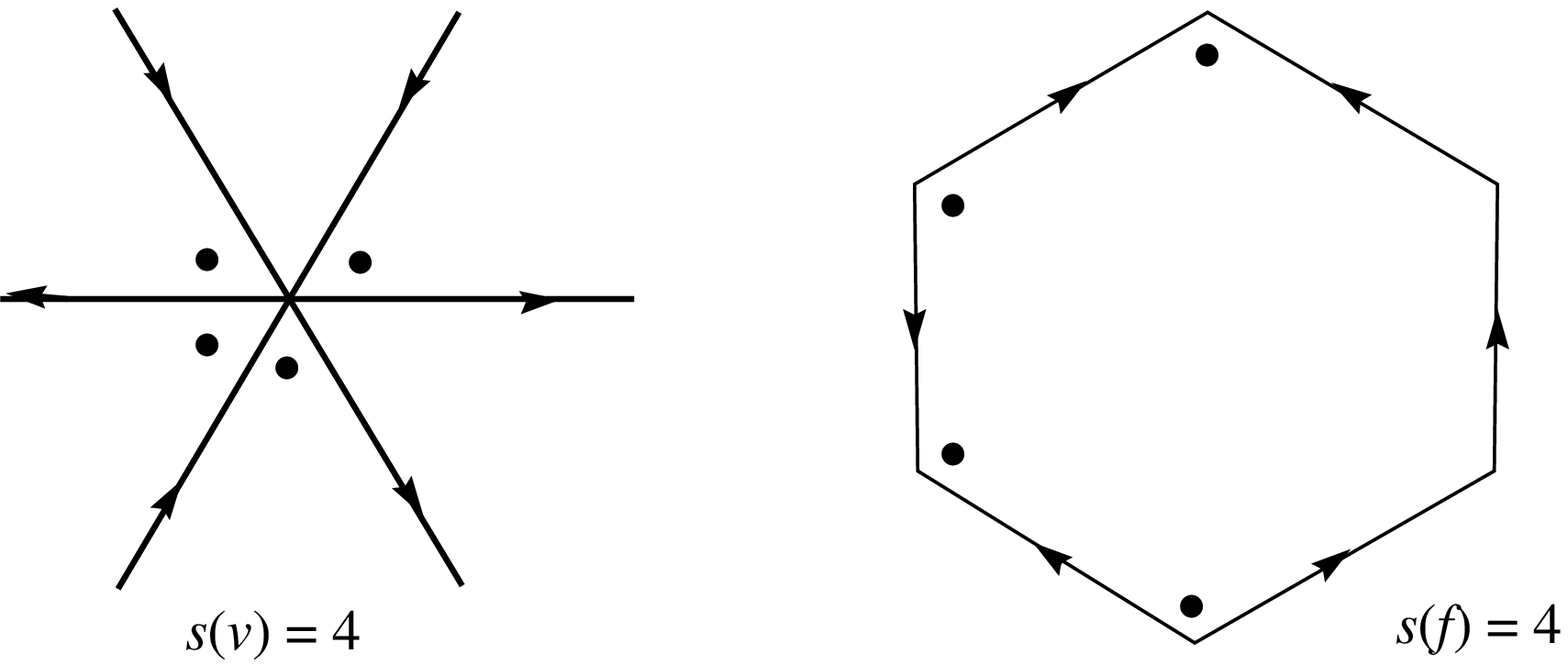}{10}

Define the {\it index\/} of a vertex or face by 
$$I(v) = 1- {s(v)\over2}\quad ,\quad I(f) = \chix(f)- {s(f)\over2}\ .$$
In particular, a vertex of positive index is a sink or source 
$(I(v)=1)$, and a face of positive index is a disk whose 
boundary is a cycle $(I(f)=1)$. 

\proclaim{Lemma 5.3} 
{\rm (Glass).} 
Let $\Omega$ be a directed graph in the closed surface $S$. 
Then $\sum I(v) + \sum I(f) = \chix (S)$.
\endproclaim 

\demo{Proof} 
(See \cite{Gl}.) 
Each corner between adjacent edges at a vertex contributes exactly 1 to 
$\sum s(v) + \sum s(f)$. 
Hence 
$$\eqalignno{
\chix (S) & = V-E+{\textstyle\sum}\chix(f)\cr 
&= V- {\#\text{ corners}\over2} + {\textstyle\sum}\chix(f)\cr
&= V - {(\sum s(v) + \sum s(f))\over2} +{\textstyle\sum}\chix(f) \cr 
&= {\textstyle\sum}I(v) + {\textstyle\sum}I(f)\ .&\qed\cr}$$
\enddemo

A label $i$ which at a fat vertex lies immediately between oppositely oriented 
edges of $\Gamma (\tau)$ is a {\it switch label\/} of $\tau$. 
(Equivalently, $i$ is a switch label if and only if 
$\tau |(i-1,i)$ and $\tau|(i,i+1)$ 
have opposite signs.) 
If the edges of $\Gamma(\tau)$ adjacent to $i$ are oriented in a clockwise 
(resp. anticlockwise) direction around $i$, then $i$ is a {\it clockwise\/} 
(resp. {\it anticlockwise\/}) switch label. 
Note that this is well-defined, i.e., independent of the sign of the fat 
vertex under consideration; see Figure~9. 

Let $C(\tau),A(\tau)\subset \bq$ be the set of clockwise and anticlockwise 
switch labels of $\tau$, respectively. 

As an example, for the type $\tau$ illustrated in Figure~9, 
$C(\tau) = \{2,4\}$, $A(\tau)=\{1,3\}$. 

An edge $e$ of $G_P$ is a {\it clockwise\/} (resp. {\it anticlockwise\/}) 
{\it switch edge\/} if the labels at the endpoints of $e$ both belong to 
$C(\tau)$ (resp. $A(\tau)$). 
A {\it switch edge\/} is an edge of $G_P$ that is either a clockwise or an 
anticlockwise switch edge. 

Note that the faces of $\Gamma (\tau)$ are in one-one correspondence with 
the edges of $G_P$, and that under this correspondence, faces of 
$\Gamma(\tau)$  of index~1 correspond to switch edges of $G_P$. 
See Figure~11. 

\bigskip

\Fig{.75}{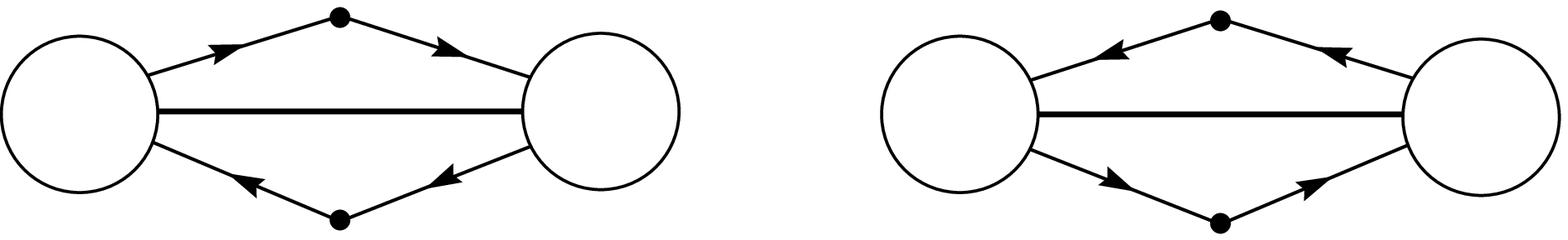}{11}

Now let $\tau$ be a non-trival $\bq$-type, and write $c(\tau) = |C(\tau)|$, 
$a(\tau) = |A(\tau)|$. 
Thus $c(\tau) = a(\tau)>0$. 

If $v$ is a fat vertex of $\Gamma (\tau)$, then 
$$s(v) = \Delta (c(\tau) + a(\tau))= 2\Delta c(\tau)\ .$$ 
Therefore 
$$I(v) = 1 - {s(v)\over2} = 1-\Delta c(\tau)\ .$$

Assume that $G_P$ does not represent $\tau$. 
Then by Proposition~5.2, $\sum_{v\text{ dual}} I(v)\le 0$. 
Therefore, by Lemma~5.3, 
$${\textstyle\sum}I(f)\ge\chix (\hatP) - {\textstyle\sum}_{v\text{ fat}} I(v) = 
\chix (\hatP) - p(1-\Delta c(\tau))\ .$$ 

Note that since $c(\tau) \ge1$ and $\Delta >1-\chix (\hatP)/p$ by hypothesis, 
the last quantity above is $>0$. 

By the remark above following the definition of a switch edge, 
the number of switch edges of $G_P \ge \sum I(f)$. 
Hence we may assume without loss of generality that there are at least 
$\sum I(f)\over2$ clockwise switch edges, and so the number of endpoints 
of clockwise switch edges 
$\ge \sum I(f) \ge p\Delta c(\tau) - (p-\chix (\hatP))$. 

Now there are $\Delta c(\tau)$ clockwise switch labels at each vertex of 
$G_P$, giving a total of $p\Delta c(\tau)$ such labels in $G_P$, and so 
there are at most $(p-\chix (\hatP))$ of these labels that are not endpoints 
of clockwise switch edges of $G_P$. 
Let $\Lambda$ be the subgraph of $G_Q$ consisting of those edges that 
correspond to clockwise switch edges of $G_P$. 
(In particular, the vertices of $\Lambda$ correspond to a subset of 
$C(\tau)$.) 
As noted above, $\Lambda$ is non-empty.    
Also, there are at most $(p-\chix (\hatP))$ occurrences of labels at vertices 
of $\Lambda$ that are not endpoints of edges of $\Lambda$. 
Hence if 
$$\align
\text{\it all elements of $C(\tau)$ have the same sign,}&\qquad (*c)\cr 
\text{and {\it all elements of $A(\tau)$ have the same sign,}}&\qquad (*a)
\endalign$$
then all the vertices of $\Lambda$ have the same sign, and so (a connected 
component of) $\Lambda$ is a $(p-\chix (\hatP))$-web, as desired. 
(Note that we need to make both assumptions $(*c)$ and $(*a)$, as we have 
no control over whether there are at least $\sum I(f)\over2$ clockwise or 
anticlockwise switch edges.) 

If at least one of $(*c)$ and $(*a)$ does not hold, then we take {\it derived 
types\/} until both do hold, and then work backwards. 
We now briefly discuss this inductive procedure. 

For this, we have to work in the following more general setting. 
Consider a non-empty set of labels $L\subset \bq$. 
Associated with $L$ is the set of {\it $L$-intervals\/}, which are the 
intervals $(\ell_1,\ell_2)$ between successive elements of $L$. 
An {\it $L$-type\/} is defined in the obvious way, to be an $|L|$-type where 
the $i$th co-ordinate is formally associated with the $i$th $L$-interval. 

Define $G(L)$ to be the subgraph of $G_P$ consisting of all edges of $G_P$ 
with at least one endpoint label belonging to $L$. 
Thus $G(L)$ and $G_P$ have the same vertices, and $G(\bq) =G_P$. 

Note that although the corners of the faces of $G(L)$ are not necessarily 
$L$-intervals, each is contained in a unique $L$-interval. 
Hence, for each disk face $D$ of $G(L)$, by taking the algebraic sum of the 
$L$-intervals containing the corners of $D$ we get an element 
$[D]\in\zed^{|L|}$. 
Then, very much as before, we say that $D$ {\it represents\/} the $L$-type 
$\tau$ if and only if 
\roster
\item"(1)" $[D]$ represents $\tau$ (in the sense of Section 3); and 
\item"(2)" ({\it coherence condition\/}) no $L$-interval occurs as the 
$L$-interval containing a corner of $D$ at two vertices of $D$ of opposite 
sign. 
\endroster

We define the directed graph $\Gamma (\tau)$ as before (starting with 
$G(L)$), orienting each edge according to the restriction of $\tau$ to the 
$L$-interval containing the corresponding corner (and the sign of the fat 
vertex in question). 
Proposition~5.2 continues to hold with $G_P$ replaced by $G(L)$. 

Let $\tau$ be a non-trivial $L$-type. 
We have $C(\tau), A(\tau)\subset L$, the set of clockwise (resp. 
anticlockwise) switches of $\tau$, as before. 
The argument we have just given above carries over verbatim to our present 
more general setting, so we have the following lemma. 

\proclaim{Lemma 5.4} 
Let $\tau$ be a non-trivial $L$-type, satisfying $(*c)$ and $(*a)$, such 
that $G(L)$ does not represent $\tau$. 
Then $G_Q$ contains a $(p-\chix (\hatP))$-web.
\endproclaim 

The {\it derivative\/} of a  non-trivial $L$-type $\tau$ is the 
$C(\tau)$-type $d\tau$ defined as follows. 
Let $(c,c')$ be a $C(\tau)$-interval. 
Then $(c,c')$ contains a unique element $a\in A(\tau)$. 
We define $d\tau$ by 
$$d\tau | (c,c') = \text{sign }a\ .$$

\example{Example} 
Taking $q=16$ and starting with the $\bq$-type $\tau = + - + + - + - - + - 
+++ - - +$, Figure~12 below illustrates how one obtains $d\tau$ and 
$d^2\tau$. 
\endexample

\topinsert
{\epsfysize=7.5truein 
\centerline{\epsfbox{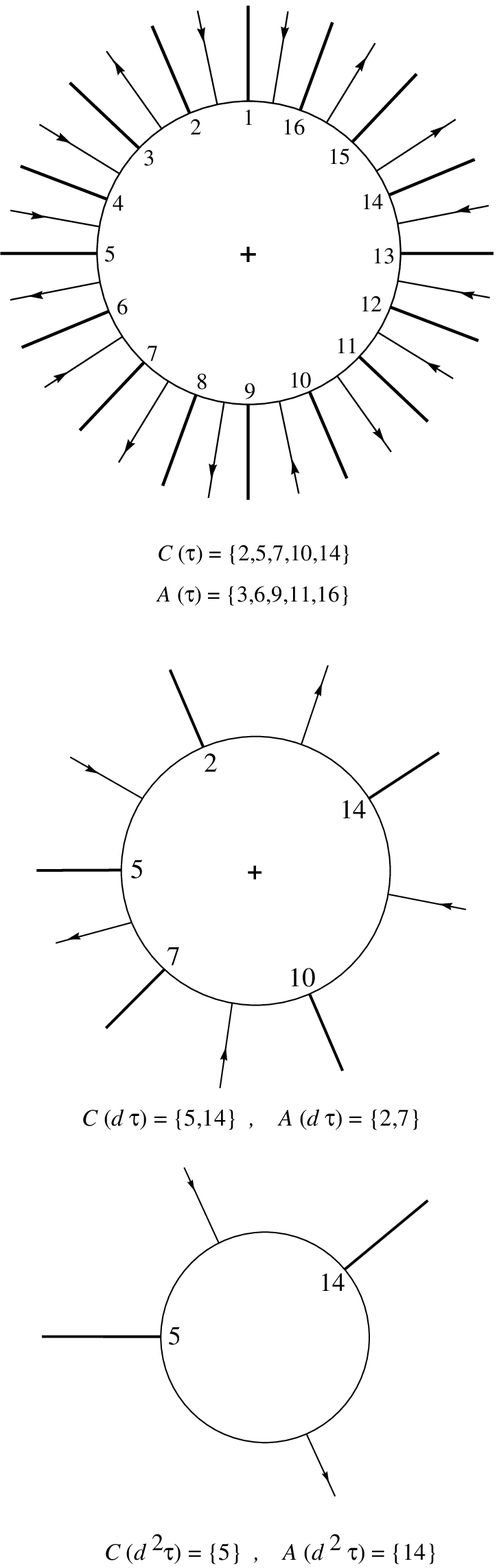}}}
\centerline{\bf Figure 12}
\endinsert

The key property of the definition of the derivative of a type is the 
following proposition, which says that to find a representative of $\tau$ 
it is enough to find a representative of $d\tau$. 

\proclaim{Proposition 5.5} 
Let $\tau$ be a non-trivial $L$-type. 
Then any face of $G(C(\tau))$ representing $d\tau$ contains a face of 
$G(L)$ representing $\tau$.
\endproclaim 

\demo{Sketch Proof} 
We shall give enough of the proof to at least explain the reason for the 
definition  of $d\tau$. 

Let $D$ be a face of $G(C(\tau))$ representing $d\tau$. 
Let $\Gamma_D (\tau) = \Gamma (\tau) \cap D$.
\enddemo 

\proclaim{Claim} 
$\Gamma_D(\tau)$ has either a face of index $1$ or a dual vertex of index $1$.
\endproclaim 

\demo{Proof of Claim} 
This is proved by doubling $(D,\Gamma_D(\tau))$ along the boundary, and 
applying Lemma~5.3 to the double $2\Gamma_D(\tau)\subset 2D\cong S^2$. 
Note that the faces of $2\Gamma_D(\tau)$ are the faces of $\Gamma_D(\tau)$ 
(each appearing twice), together with a face $f_e$ for each edge 
$e\subset \partial D$ (with $e\subset f_e$). 
Hence the conclusion will follow if we show that 
$$|\{ \text{fat vertices of $2\Gamma_D(\tau)$ of index 1}\}| 
\le |\{ \text{faces $f_e$ of index $-1$}\}|\ .$$
But this in turn follows from the fact that, by definition of $G(C(\tau))$, 
each $e\subset \partial D$ has at least one endpoint label $\in C(\tau)$ 
(see \cite{GL2, p.398} for details). 

A dual vertex of $\Gamma_D(\tau)$ of index 1 will correspond to a face 
of $G(L)$, contained in $D$, representing $\tau$. 
So, by the claim above, we are done once we show that no face of 
$\Gamma_D(\tau)$ has index~1. 
Such a face corresponds to a switch edge $e$ of $G(L)$. 
Now $e$ cannot be a clockwise switch edge, for such edges are edges of 
$G(C(\tau))$, and hence cannot lie in the interior of the face $D$ of 
$G(C(\tau))$. 
So assume that $e$ is an anticlockwise switch edge, with 
endpoint labels $a_1,a_2\in A(\tau)$ at vertices $v_1,v_2$ respectively; 
see Figure~13. 
Then, by the parity rule, $(\text{sign }v_1)(\text{sign }a_1) \ne 
(\text{sign }v_2)(\text{sign }a_2)$. 
But, recalling the definition of $d\tau$, this contradicts the fact 
that $D$ represents $d\tau$.\qed
\enddemo 

\bigskip

\Fig{.75}{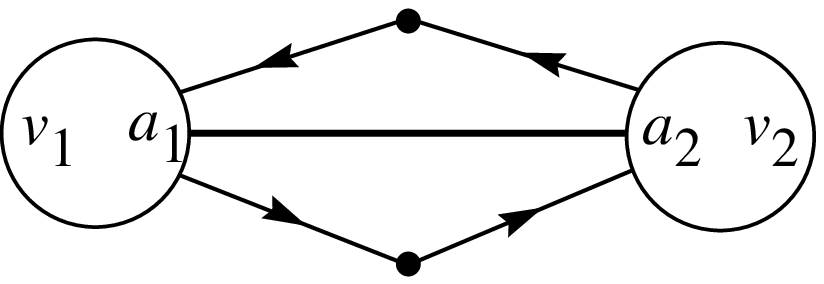}{13}

Returning to the proof of Theorem 5.1, suppose that one of the conditions 
$(*c),(*a)$ fails for the $\bq$-type $\tau$. 
Note that $C(\tau) = A(-\tau)$. 
Hence, if we define 
$$\tau_0 = \cases 
\hphantom{-}\tau\ ,&\text{if $\tau$ does not satisfy $(*a)$}\cr 
-\tau\ ,&\text{if $\tau$ satisfies $(*a)$ but not $(*c)$,}
\endcases$$ 
then $(*a)$ fails for $\tau_0$. 
Hence $d\tau_0$ is a non-trivial type. 
If $d\tau_0$ satisfies $(*c)$ and $(*a)$, stop. 
If not, define 
$$\tau_1 = \cases 
\hphantom{-}d\tau_0 \ ,&\text{if $d\tau_0$ does not satisfy $(*a)$}\cr 
-d\tau_0\ ,&\text{if $d\tau_0$ satisfies $(*a)$ but not $(*c)$}
\endcases$$ 
Continuing in this way, we obtain a sequence $\tau_0,\tau_1,\ldots,\tau_n$, 
$n\ge1$, where $\tau_0 =\pm\tau$, $\tau_i = \pm d\tau_{i-1}$ is a non-trivial 
$C(\tau_{i-1})$-type, $1\le i\le n$, and $\tau_n$ satisfies $(*c)$ and $(*a)$.

If $G(C(\tau_{n-1}))$ does not represent $\tau_n$, then, by Lemma~5.4, 
$G_Q$ contains a $(p-\chix (\hatP))$-web. 

If $G(C(\tau_{n-1}))$ does represent $\tau_n$, then successive applications of 
Proposition~5.5 show that $G_P$ represents $\tau$. 

This completes our sketch of the proof of Theorem 5.1. 
\medskip

We have seen, in Section 4, that if $G_P$ represents all types then we get 
useful topological information about $(M(\beta),\hatQ)$. 
What if $G_Q$ contains a $(p-\chix (\hatP))$-web? 
One idea is that a web might give rise to Scharlemann cycles in $G_Q$ 
(which then in turn give useful topological information about 
$(M(\alpha),\hatP))$. 
In order to make this work, however, we must impose a further condition. 
Namely, we define a {\it great\/} $k$-web in $G_Q$ to be a $k$-web $\Lambda$ 
which is contained in the interior of a disk $D_\Lambda \subset\hatQ$ with 
the property that any vertex of $G_Q$ lying in $D_\Lambda$ is a vertex 
of $\Lambda$. 
See Figure~14.

\bigskip

\Fig{3.0}{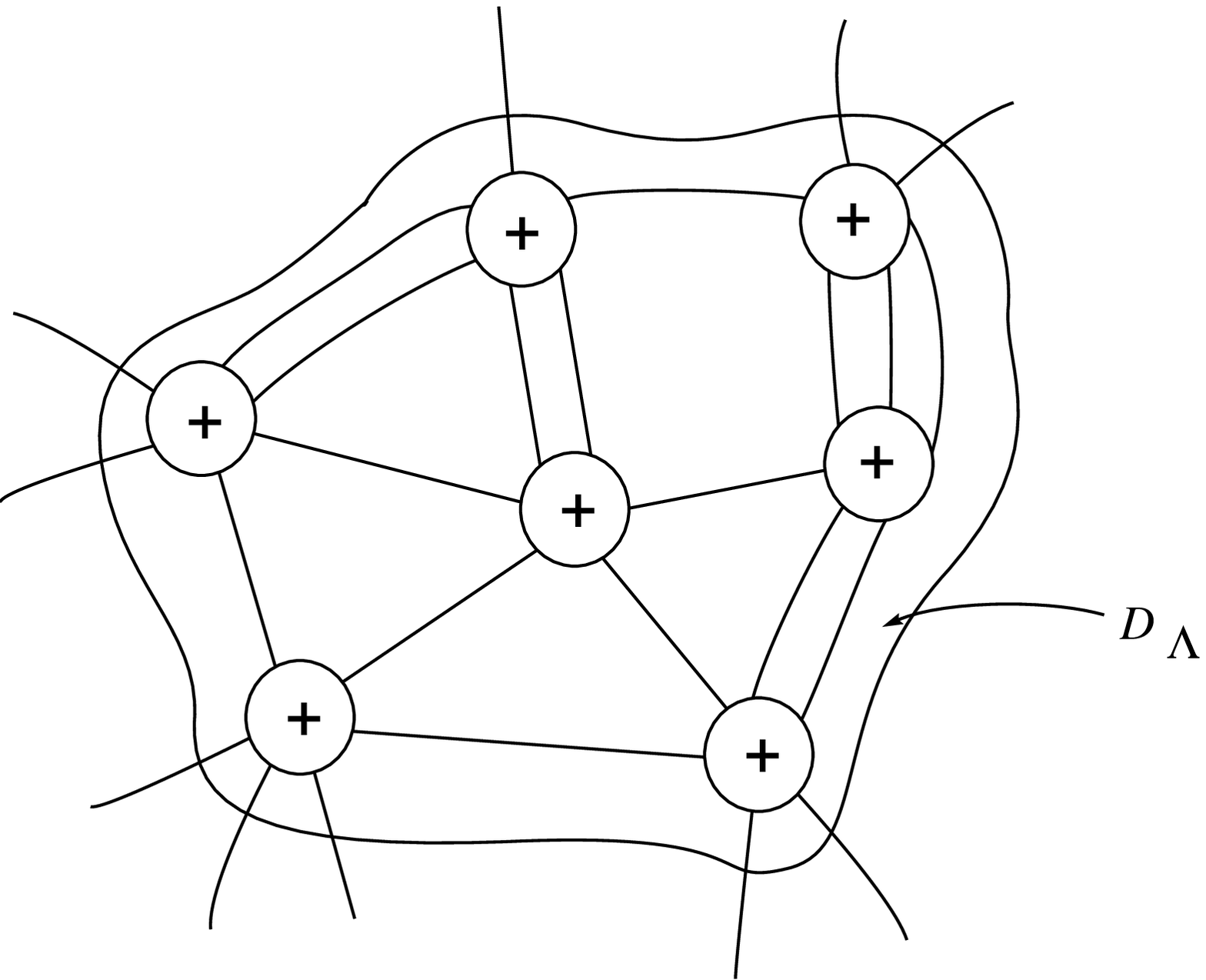}{14}
\noindent
To guarantee that such objects exist, we need to specialize to the case 
where $\hatQ$ is a 2-sphere. 
It turns out that one can then prove the following theorem. 

\proclaim{Theorem 5.6} 
Assume that $\Delta >1-\chix (\hatP)/p$, and that $\hatQ\cong S^2$. 
Then either $G_P$ represents all types or $G_Q$ contains a great 
$(p-\chix (\hatP))$-web.
\endproclaim 

The proof proceeds as follows. 
By Theorem 5.1, if $G_P$ does not represent all types then $G_Q$ contains 
a $(p-\chix(\hatP))$-web $\Lambda$. 
Let $U$ be a component 
of $\hatQ-\nhd (\Lambda)$, and let $D$ be the disk $\hatQ-U$, containing 
$\Lambda$. 
If $\Lambda$ is not a great web, then there are vertices of $G_Q$ in $D$ that 
are not vertices of $\Lambda$. 
These vertices correspond to a non-empty subset $L\subset \bq$. 
One then proves by induction on $|L|$ that either $G_Q$ contains a great 
$(p-\chix (\hatP))$-web or $G(L)$ represents all $L$-types. 
(The last step of the induction, which is the desired statement, is 
the degenerate case where $\Lambda =\emptyset$, $L=\bq$, and $D$ is a disk 
in $\hatQ$ containing $G_Q$.) 
This inductive argument is quite subtle, and in particular involves 
relativizing the notion of the derivative of a type. 
We refer to  \cite{GL3, proof of Theorem~2.5} (which does the case where 
$\hatP$ is a torus), and thence to \cite{GL2, Section~2}, for more details. 

\head 6. Applications\endhead 

In this section we sketch the proofs of some results which use the methods 
we have discussed. 
The first three are about Dehn surgery on knots in $S^3$, and the proofs 
of these will all be based on the following statement, which is an 
immediate consequence of Theorem~5.6 and (Ta) of Section~4. 

\proclaim{Theorem 6.1} 
Suppose that $M(\beta)\cong S^3$ and that $\hatQ$ is a Heegaard 2-sphere. 
If $\Delta >1-\chix (\hatP)/p$ then $G_Q$ contains a great 
$(p-\chix(\hatP))$-web. 
\endproclaim 

When $M(\beta) \cong S^3$, we shall make a slight change of notation, writing 
$K$ for $K_\beta$, a non-trivial knot in $S^3$, $M_K$ for $M$, the exterior 
of $K$, and $\mu_K$ for $\beta$, the meridian of $K$. 

The following theorem is proved in \cite{GL2}. 
It implies that knots in $S^3$ are determined by their complements. 

\proclaim{Theorem 6.2} 
Let $K$ be a non-trivial knot in $S^3$. 
If $M_K(\alpha) \cong S^3$ then $\alpha = \mu_K$. 
\endproclaim 

\demo{Proof} 
Suppose that $M_K (\alpha) \cong S^3$, with $\Delta (\alpha,\mu_K)\ge1$. 
Let $\hatP$ be a Heegaard 2-sphere in $M_K(\alpha)$. 

By Theorem 6.1, $G_Q$ contains a great $(p-2)$-web, $\Lambda$. 
Hence there exists a label (in fact at least two labels) $i\in \bp$ such that 
every occurrence of $i$ at a vertex of $\Lambda$ is the endpoint of an edge 
of $\Lambda$. 
Moreover, since all the vertices of $\Lambda$ have the same sign, no edge of 
$\Lambda$ can have the same label at both endpoints, by the parity rule. 
Hence, by always leaving a vertex of $\Lambda$ along an edge of $\Lambda$ 
with label $i$ at that vertex, we see that $\Lambda$ contains an 
{\it $i$-cycle\/}, that is, a cycle of edges which can be consistently 
oriented  so that the tail of each edge has label $i$ (and all the vertices 
in the cycle have the same sign). 
See Figure~15. 

\bigskip

\Fig{2.0}{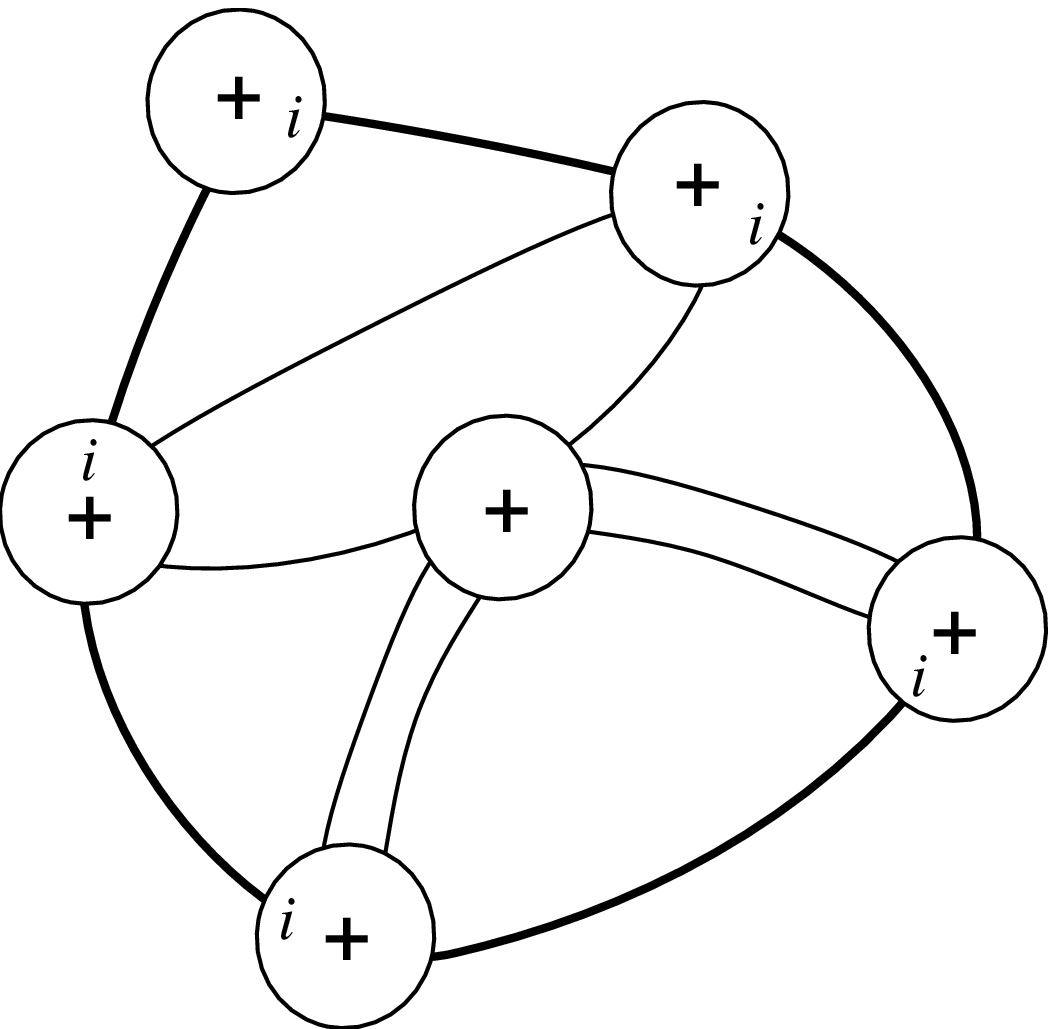}{15}
\noindent
Now it is easy to see that a cycle of this kind which is innermost   in 
the disk $D_\Lambda$ (allowing the label $i$ to be an arbitrary element 
of $\bp$) is a Scharlemann cycle in $G_Q$. 
Since $\hatP$ is a 2-sphere the edges of this Scharlemann cycle necessarily 
lie in a disk in $\hatP$. 
But since $\hatP$ is a Heegaard 2-sphere, this is a contradiction, by 
(Sa) in Section~4 (here, the roles of $P$ and $Q$ are interchanged).\qed
\enddemo 

\remark{Remark} 
We observed in the above proof that a great $(p-2)$-web in $G_Q$ always 
contains a Scharlemann cycle. 
So we may state
\endremark

\proclaim{Theorem 6.3} 
If $\hatP$ and $\hatQ$ are 2-spheres, and $\Delta\ge1$, then either 
$G_P$ represents all types or $G_Q$ contains a Scharlemann cycle.
\endproclaim 

If $\Delta\ge2$ then the following stronger conclusion holds. 

\proclaim{Theorem 6.4} 
If $\hatP$ and $\hatQ$ are 2-spheres, and $\Delta\ge2$, then either  
$G_P$ or $G_Q$ contains a Scharlemann cycle.
\endproclaim 

This is an immediate consequence of   \cite{CGLS, Proposition 2.5.6} 
(which deals with the case where $\hatP$ and $\hatQ$ are disks). 

Theorem 6.4 already shows that if $M(\alpha)\cong M(\beta)\cong S^3$ 
then $\Delta (\alpha,\beta)\le1$, implying that there are at most two 
inequivalent knots with homeomorphic complements (a fact first proved, 
by a different argument, in \cite{CGLS}). 
However, there are examples of graphs $G_P,G_Q$ with $\Delta =1$ for which 
the conclusion of Theorem~6.4 fails, illustrating the necessity of 
introducing the more general concept of representing all types. 

The next application addresses the question  of when a reducible manifold 
can be obtained by Dehn surgery on a knot in $S^3$. 

\proclaim{Theorem 6.5} 
Let $K$ be a  non-trivial knot in $S^3$. 
If $M_K(\alpha)$ is reducible then 
\roster
\item"(1)"    $M_K(\alpha)$ has a lens space summand; and 
\item"(2)" $\Delta (\alpha,\mu_K)=1$.
\endroster
\endproclaim 

Part (1) was proved in \cite{GL2}; part (2) was proved in \cite{GL1} (by an 
argument different from that outlined below). 

Part (1) implies the Property $R$ Conjecture, proved by Gabai \cite{Ga2}. 

\proclaim{Corollary 6.6 {\rm (Gabai)}} 
If $K$ is a non-trivial knot in $S^3$ then $M_K(0)$ is irreducible.
\endproclaim 

Part (1) also implies 

\proclaim{Corollary 6.7} 
Any homology 3-sphere obtained by Dehn surgery on a knot in $S^3$ is prime.
\endproclaim 

We remark that Auckly, using results from gauge 
theory, has shown that there are prime homology 
3-spheres which cannot be obtained by Dehn surgery on any knot in 
$S^3$ \cite{A}. 

The situation described in Theorem 6.5 can occur: 
if $K$ is a torus knot or cable knot then there is a slope $\alpha$ such 
that $M_K(\alpha)$ is reducible. 
However, the Cabling Conjecture \cite{GS} asserts that these are the only 
knots with this property. 
This is still open. 

\demo{Sketch Proof of Theorem 6.5} 
Here we take $\hatP$ to be an essential 2-sphere in $M_K(\alpha)$. 
(Note that $M_K$ is irreducible.) 

By Theorem 6.1, $G_Q$ contains a great $(p-2)$-web $\Lambda$. 
As in the proof of Theorem~6.2 above, this implies that $G_Q$ contains a 
Scharlemann cycle, and hence that $M_K(\alpha)$ has a lens space summand 
(see (Sc) in Section~4). 
This proves (1). 

To prove (2), assume that $\Delta (\alpha,\mu_K)>1$. 
Then, for homological reasons, $\hatP$ separates $M_K(\alpha)$, and 
hence $p$ is even. 

If $p=2$ then $P$ is an annulus, and hence $K$ is a  cable of a knot $K'$ 
(which may be trivial). 
Now it is not hard to show (see the first four lines of the proof of 
Theorem~3 in \cite{GL1}) that if $M_K(\alpha)$ is reducible for some $\alpha$ 
with $\Delta (\alpha,\mu_K)>1$ then $M_{K'}(\alpha')$ is reducible for 
some $\alpha'$ with $\Delta (\alpha',\mu_{K'})>1$. 
Hence, by induction it suffices to prove the result when $K$ is not cabled. 

So assume $p\ge4$. 
As in the proof of Theorem 6.2 above, there is a label in $\bp$ such 
that each of the $\Delta$ occurrences of that label at each vertex 
of $\Lambda$ is the endpoint of an edge of $\Lambda$. 
One can then show, using the fact that $\Delta \ge2$, that $\Lambda$ contains 
Scharlemann cycles of $G_Q$ on distinct $\bp$-intervals $(i,i+1)$, 
$(j,j+1)$ 
(see \cite{GL4, Proof of Theorem 2.3}). 
This allows one to construct a new essential 2-sphere $\hatP' \subset 
M_K(\alpha)$ such that $|\hatP'\cap K_\alpha| <p$ 
(see \cite{GL4, Theorem 2.4}) 
contradicting the minimality of $p$.\qed
\enddemo 

For our last application to Dehn surgery on knots in $S^3$ we consider 
the situation where $M_K(\alpha)$ contains an incompressible torus. 
Here, we must assume that $M_K$ contains no essential torus, or equivalently, 
that $K$ is not a satellite knot. 

\proclaim{Theorem 6.8 \cite{GL3}} 
Let $K$ be a knot in $S^3$ that is not a satellite knot. 
If $M_K(\alpha)$ contains an incompressible torus then 
$\Delta (\alpha, \mu_K) \le 2$.  
\endproclaim 

Examples with $\Delta (\alpha,\mu_K)=2$ have been constructed by 
Eudave-Mu\~noz \cite{E2}. 
There are many examples with $\Delta (\alpha,\mu_K)=1$. 

\demo{Sketch Proof of Theorem 6.8} 
Let $\hatP$ be an incompressible torus in $M_K(\alpha)$, and assume (for 
a contradiction) that $\Delta (\alpha,\mu_K)\ge3$. 
In particular, this implies that $\hatP$ separates $M_K(\alpha)$, so that 
$p$ is even. 

If $p=2$ then one can show that the knot $K$ is strongly invertible 
(see \cite{GL3, Section~8}) 
in which case the theorem is proved by Eudave-Mu\~noz \cite{E1}. 

So assume that $p\ge4$. 
By Theorem~6.1, $G_Q$ contains a great $p$-web $\Lambda$. 
An easy euler characteristic argument 
(see \cite{GL3, Section~4}) shows that there are at least 
four labels $i\in \bp$ such that $\Lambda$ contains an {\it $i$-bigon\/}; 
see Figure~16. 
Now such an $i$-bigon is either a Scharlemann cycle, or contains within it 
an {\it extended Scharlemann cycle\/}, that is, a Scharlemann cycle flanked 
by bigons; 
see Figure~17. 
However, a topological argument shows that an extended Scharlemann cycle in 
$G_Q$ would give rise to a new incompressible torus $\hatP'$ in $M_K(\alpha)$ 
with $|\hatP'\cap K_\alpha| <p$, contradicting the minimality of $p$ 
(see \cite{GL3, Theorem 3.2}).   
We conclude that $\Lambda$ contains two Scharlemann cycles $D,D'$ of 
length~2 on disjoint $\bp$-intervals $(i,i+1)$ and $(i',i'+1)$. 
By (Sb) in Section~4, the edges of neither $D$ nor $D'$ lie in a disk 
in $\hatP$, hence each pair of edges forms an essential loop on $\hatP$; 
see Figure~18. 
Shrinking the 1-handle $H_{(i,i+1)}$ to its core has the effect of gluing 
together the two $(i,i+1)$ corners of $D$ 
(see Figure~19), 
giving a M\"obius band~$B$. 
Similarly, $D'$ gives rise to a M\"obius band $B'$. 
Since $\partial B$ and $\partial B'$ are parallel on $\hatP$ we may join 
$B$ and $B'$ by an annulus in $\hatP$ to get a Klein bottle $F$ 
in $M_K(\alpha)$. 
We now replace the torus $\hatP$ by the boundary $\hatP_0$ of a regular 
neighborhood of $F$. 
(Note that $\hatP_0$ is incompressible, since $M_K(\alpha)$ contains an 
incompressible torus by hypothesis, and is irreducible by Theorem~6.5(2).) 
But similar combinatorial arguments, applied to $P_0\cap Q$, in this case 
lead to a contradiction 
(see \cite{GL3, Section~6}).\qed 
\enddemo 

\bigskip

\Fig{.7}{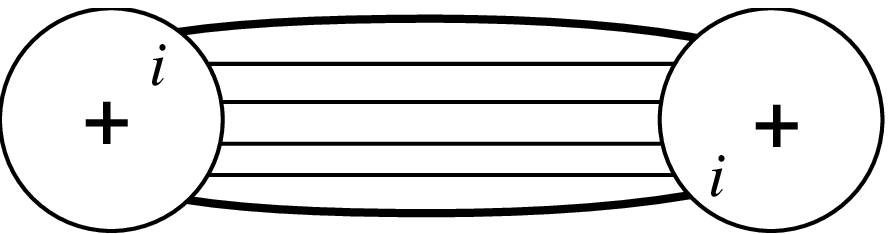}{16}

\bigskip

\Fig{.7}{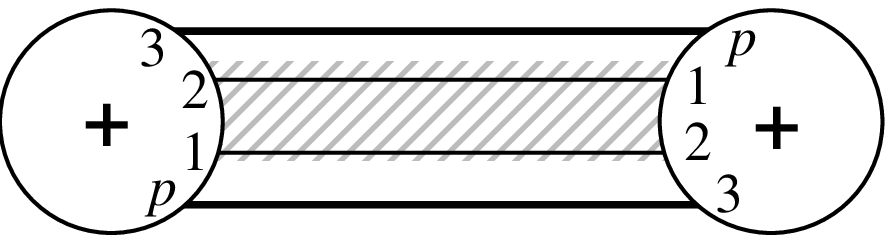}{17}

\bigskip

\Fig{3.0}{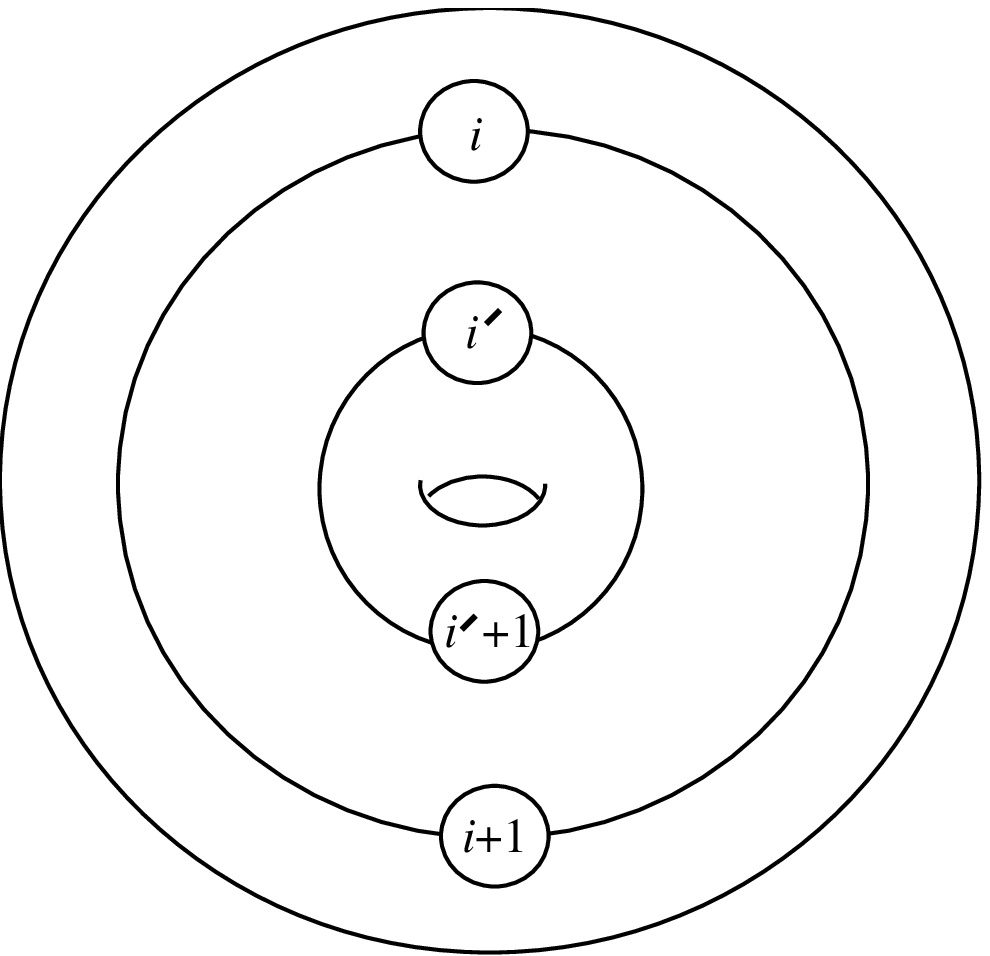}{18}

\bigskip

\Fig{.75}{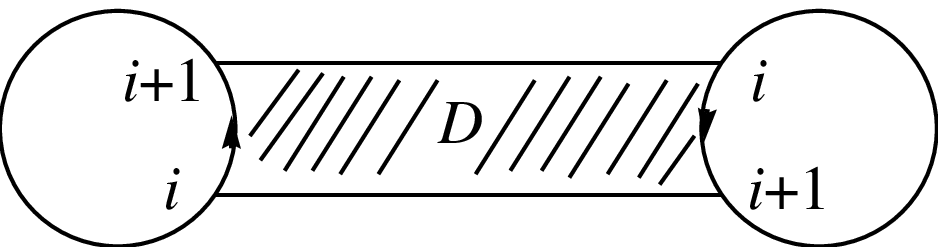}{19}

Finally, we show how the methods we have discussed give some information 
about the question of when $M(\alpha)$ and $M(\beta)$ are reducible. 
Let us say that a 3-manifold is {\it $\que$-reducible\/} if it contains 
a $\que$-essential 2-sphere. 

\proclaim{Theorem 6.9} 
\roster 
\runinitem"(i)" If $M(\alpha)$ and $M(\beta)$ are $\que$-reducible then 
$\alpha =\beta$. 
\item"(ii)" If $M(\alpha)$ is reducible and $M(\beta)$ is $\que$-reducible, 
with $\alpha\ne\beta$, 
then $M(\alpha)$ has a lens space summand.
\item"(iii)" If $M(\alpha)$ and $M(\beta)$ are reducible, 
with $\alpha\ne\beta$, then either both 
$M(\alpha)$ and $M(\beta)$ have summands $M_\alpha,M_\beta$ with 
$H_1(M_\alpha)$ and $H_1(M_\beta)$ finite and non-zero, or one of 
$M(\alpha),M(\beta)$ has a lens space summand.
\endroster
\endproclaim 

\remark{Remarks}
(1) Part (i) is a special case of a result of Scharlemann \cite{S3} proved 
using the theory of sutured manifolds. 
The case where $M(\alpha)$ and $M(\beta)$ are connected sums of copies 
of $S^1\times S^2$ was proved by Gabai \cite{Ga1}. 

(2) Scharlemann \cite{S3} proves that under the hypotheses of (ii), 
$M(\beta) \cong S^1\times S^2\, \#\, W$, and $M(\alpha) \cong W'\,\#\, W''$, 
where $W,W'$ and $W''$ are $\que$-homology spheres, and that either $W'$ 
is a lens space or $H_1(W)\ne0$. 

(3) It is shown in \cite{GL4}, by an elaboration of the methods discussed here, 
that if $M(\alpha)$ and $M(\beta)$ are reducible then $\Delta (\alpha,\beta)
\le1$. 
A different proof has been given by Boyer and Zhang \cite{BZ}. 
\endremark 

\demo{Proof of Theorem 6.9} 
(i) 
Let $\hatP,\hatQ$ be $\que$-essential 2-spheres in $M(\alpha)$, 
$M(\beta)$ respectively with $p$ and $q$ minimal. 
By Theorem~6.3, either $G_P$ represents all types or $G_Q$ contains a 
Scharlemann cycle. 
But the first contradicts the minimality of $q$ (see (Td) in Section~4), 
while the second contradicts the minimality of $p$ (see (Sd) in Section~4). 

(ii) 
Let $\hatQ$ be as in (i), and let $\hatP$ be an essential 2-sphere 
in $M(\alpha)$. 
By Theorem~6.3 (and (Td)) $G_Q$ contains a Scharlemann cycle. 
Hence $M(\alpha)$ has a lens space summand. 

(iii) 
Let $\hatP,\hatQ$ be essential 2-spheres in $M(\alpha)$, $M(\beta)$ 
respectively. 
The result then follows from Theorem~6.3, applied as stated and also 
with $P$ and $Q$ interchanged, (together with (Sc) and (Tc) in Section~4).\qed
\enddemo

\enddocument